\documentclass[11pt,letterpaper]{amsart}
\usepackage{amsmath,amsthm,amsxtra,amssymb,latexsym,mathrsfs,mathtools}
\usepackage[a4paper,top=2.5cm,bottom=2cm,left=2.5cm,right=2.5cm,marginparwidth=1.75cm]{geometry}
\usepackage{verbatim}
\usepackage{amscd}
\usepackage[all]{xy}
\usepackage{eucal}
\usepackage{hyperref}
\usepackage{color}
\usepackage{float}
\usepackage{tikz}
\usepackage{comment}
\usepackage{enumitem}

\newtheoremstyle{plain-alt}% name of the style to be used
  {}% measure of space to leave above the theorem. E.g.: 3pt
  {}% measure of space to leave below the theorem. E.g.: 3pt
  {\itshape}% name of font to use in the body of the theorem
  {}% measure of space to indent
  {}% name of head font
  {. }% punctuation between head and body
  {0em}% space after theorem head; " " = normal interword space
  {\thmnumber{#2}.\,\,\thmname{{\it #1}}\thmnote{ (#3)}}% Manually specify head
\newtheoremstyle{definition-alt}% name of the style to be used
  {}% measure of space to leave above the theorem. E.g.: 3pt
  {}% measure of space to leave below the theorem. E.g.: 3pt
  {\normalfont}% name of font to use in the body of the theorem
  {}% measure of space to indent
  {}% name of head font
  {. }% punctuation between head and body
  {0em}% space after theorem head; " " = normal interword space
  {\thmnumber{#2}.\,\,\thmname{{\it #1}}\thmnote{ (#3)}}% Manually specify head
\theoremstyle{plain-alt}

\theoremstyle{plain}
\swapnumbers
\newtheorem{lemma}{Lemma}[section]
\newtheorem{prop}[lemma]{Proposition}
\newtheorem{corollary}[lemma]{Corollary}
\newtheorem{theorem}[lemma]{Theorem}

\newtheorem*{thm}{Theorem}

\theoremstyle{definition}
\newtheorem{definition}[lemma]{Definition}
\newtheorem{construction}[lemma]{Construction}
\newtheorem{remark}[lemma]{Remark}

\newtheorem{example}[lemma]{Example}

\begin{document}

\newcommand{\wt}{\operatorname{wt}}
\newcommand{\norm}{\operatorname{N}}
\newcommand{\tr}{\operatorname{tr}}
\newcommand{\Lie}{\operatorname{Lie}}
\newcommand{\Aut}{\operatorname{Aut}}
\newcommand{\Out}{\operatorname{Out}}
\newcommand{\Der}{\operatorname{Der}}
\newcommand{\End}{\operatorname{End}}
\newcommand{\Int}{\operatorname{Int}}
\newcommand{\Ad}{\operatorname{Ad}}
\newcommand{\ad}{\operatorname{ad}}
\newcommand{\GL}{\operatorname{GL}}
\newcommand{\SL}{\operatorname{SL}}
\newcommand{\ind}{\operatorname{ind}}
\newcommand{\res}{\operatorname{res}}
\newcommand{\differential}{\operatorname{d}}
\newcommand{\Supp}{\operatorname{Supp}}
\newcommand{\Iso}{\operatorname{Iso}}
\newcommand{\Exp}{\operatorname{Exp}}
\newcommand{\Aff}{\operatorname{Aff}}
\newcommand{\Pic}{\operatorname{Pic}}
\newcommand{\SO}{\operatorname{SO}}
\newcommand{\U}{\operatorname{U}}
\newcommand{\Sp}{\operatorname{Sp}}
\newcommand{\Spin}{\operatorname{Spin}}
\newcommand{\SU}{\operatorname{SU}}
\newcommand{\Ind}{\operatorname{Ind}}
\newcommand{\Hom}{\operatorname{Hom}}
\newcommand{\Nm}{\operatorname{N}}
\newcommand{\PSL}{\operatorname{PSL}}
\newcommand{\diff}{\, d}
\newcommand{\Res}{\operatorname{Res}}
\newcommand{\pr}{\operatorname{pr}}
\newcommand{\I}{\operatorname{I}}
\newcommand{\sHom}{\mathcal{H}om}
\newcommand{\sEnd}{\mathcal{E}nd}
\newcommand{\DiffEnd}{\mathcal{D}iff\text{-}\mathcal{E}nd}
\newcommand{\gr}{\operatorname{gr}}
\newcommand{\Gr}{\operatorname{Gr}}
\newcommand{\Tor}{\operatorname{Tor}}
\newcommand{\Spec}{\operatorname{Spec}}
\newcommand{\Proj}{\operatorname{Proj}}
\newcommand{\Spf}{\operatorname{Spf}}
\newcommand{\Gal}{\operatorname{Gal}}
\newcommand{\Art}{\operatorname{Art}}
\newcommand{\Frob}{\operatorname{Frob}}
\newcommand{\diag}{\operatorname{diag}}
\newcommand{\upperhalf}{\mathfrak h}
\newcommand{\Ext}{\operatorname{Ext}}
\newcommand{\Def}{\operatorname{Def}}
\newcommand{\Rep}{\operatorname{Rep}}
\newcommand{\Sym}{\operatorname{Sym}}
\newcommand{\GSp}{\operatorname{GSp}}
\newcommand{\vol}{\operatorname{vol}}
\newcommand{\depth}{\operatorname{depth}}
\newcommand{\codim}{\operatorname{codim}}
\newcommand{\Ass}{\operatorname{Ass}}
\newcommand{\cone}{\operatorname{Cone}}
\newcommand{\proj}{\operatorname{Proj}}
\newcommand{\Cl}{\operatorname{Cl}}

\newcommand{\init}{\operatorname{in}}
\newcommand{\im}{\operatorname{im}}
\newcommand{\trop}{\mathcal T}
\newcommand{\GF}{\operatorname{GF}}

\newcommand\takuya[1]{{\color{blue} \sf $\infty$ Takuya: [#1]}}
\newcommand\lara[1]{{\color{teal} \sf $\clubsuit$ Lara: [#1]}}

\sloppy

\title{Maps to toric varieties, toric degenerations and integrable systems à la Harada--Kaveh}
\author{Takuya Murata and Lara Bossinger}

\address{School of Mathematics, Institute for Research in Fundamental Sciences (IPM), Teheran, Iran}
\email{takusi@gmail.com}

\address{
Instituto de Matem\'aticas Unidad Oaxaca, 
Universidad Nacional Aut\'onoma de M\'exico,
Le\'on 2, altos, 
Centro Hist\'orico,
68000 Oaxaca,
Mexico}
\email{lara@im.unam.mx}

\date{\today}

\begin{abstract} Given a toric degeneration (a degeneration to a toric variety), over the complex numbers, we construct a surjective continuous map from a general fiber to the special fiber of the degeneration in the classical topology. The construction is a variant of one due to Goresky and MacPherson based on the Thom--Mather theory of stratified spaces.

As an application, we recover and extend the construction of integrable systems à la Harada--Kaveh in \cite{HK15}. Compared to their result, our map is constructed more explicitly and we also construct the integrable systems on the boundary strata.

This paper is a part of the authors' research on maps to toric degenerations; we refer the readers to \cite{algebraic_tex} and \cite{multi_proj} for more algebraic approaches.
\end{abstract}

\maketitle

\bigskip

\setcounter{section}{-1}
\setcounter{tocdepth}{1}

\tableofcontents

\section{Introduction}

In this paper, given a toric degeneration (see below for the definition) $\mathfrak{X}$ from a variety $X$ to a not-necessarily-normal toric variety $W$ over complex numbers, we construct a surjective map
$$\varphi : X \to W$$
that is continuous and has some more properties noted below.

More precisely, we assume $\mathfrak{X}$ is a closed subvariety of a complex smooth variety $M$ (e.g., $M = \mathbb{P}^n \times \mathbb{A}^1$) and then we construct a neighborhood $U$ of $W$ in $M$ that is an appropriate generalization of a tubular neighborhood for $W \subset M$ and then $\varphi$ is the restriction of the projection map
$$\pi : U \to W$$
to a general fiber $X$ in $\mathfrak{X}$.

A map similar to $\varphi$ was originally constructed by Goresky and MacPherson (\cite{MacPhersonGoresky} and \cite{Goresky}). 
Their construction as well as ours are based on the Thom–Mather theory of stratified spaces, which, given a stratified space (a space with stratification), constructs a system of fibrations parametrized by the strata. Roughly speaking, the above $\pi$ is obtained by composing tubular neighborhood projections for the strata on $W$. In particular, $\pi$ is a (strong) deformation retract of $W \hookrightarrow X$.

%The theory is then applied to the following setup: we have a (not-necessarily-normal) toric variety $W$ embedded in some ambient space over complex numbers (we use the classical topology); $W$ is naturally stratified by orbits. For each orbit $A$, there is a fibration $\pi_A$ from a neighborhood of $A$ to $A$ (explicitly, the neighborhood is a tubular neighborhood of $A$). Combining those $\pi_A$, we get a retraction $\pi$ from a neighborhood of $W$ to $W$. Finally, restricting $\pi$ to a general fiber of a degeneration, we get a map from a general fiber to the special fiber $W$. This is a purely topological construction. But here we use a real-algebraic-version of the theory so the resulting map $\pi$ actually reflects some geometry too.

The above map $\varphi$ and similar maps have the applications, including but not limited to:
\begin{enumerate}[label=(\arabic*)]
\item The construction of an integrable system à la Harada and Kaveh \cite{HK15}.
\item Computation of an intersection homology; thus of a perverse sheaf.
\end{enumerate}

The present paper gives the application (1), see \S\ref{sec:symplectic considerations} below. 
The application (2) is the original application of Goresky and MacPherson and is also related to the work of Clemens \cite{Clemens}, see \S\ref{sec:further} below. 
Our paper does not include the application (2), although it is related to the paper the first author is working on right now \cite{AMurata}.

\subsection{Toric degenerations} Here, we first define a toric degeneration. For a more extended review, in particular Gröbner aspects, we refer the readers to \cite{algebraic_tex}.

First let us define a toric degeneration. A degeneration is often constructed or characterized in terms of generalized extended Rees algebras (which geometrically amount to open loci in generalized blow-ups). Abstracting this, we can thus give the following definition: given an algebraic variety $X$, a {\it (flat) abstract degeneration of $X$} is a flat surjective morphism of schemes $\xi: \mathfrak{X} \to \mathbb{A}^1$ together with an isomorphism over $\mathbb{A}^1 - 0$:
$$\xymatrix{\mathfrak{X} - \xi^{-1}(0) \ar[rd]_\xi \ar[rr]^*[@]{\sim} & & X \times (\mathbb{A}^1 - 0) \ar[ld]^{\mathrm{pr}_2} \\
& \mathbb{A}^1 - 0.}\\$$
(Here $\mathbb{A}^1$ is a common choice but can be more general; e.g., \S \ref{sec:Rees algebra}) The scheme $\mathfrak{X}$ is called {\it  the total space of the degeneration} and $\xi^{-1}(0)$ called the {\it  special fiber}. Intuitively speaking, we think of a degeneration as general fibers approaching the special fiber, and the special fiber is often also called the \emph{limit} of the degeneration.

Finally, an {\it (abstract) toric degeneration} is a degeneration whose special fiber is a (not-necessarily-normal) toric variety. 

We stress that the special fiber of a toric degeneration is not a toric variety in the usual sense; in particular, it is not defined by a fan of cones. There are several approaches to not-necessarily-normal toric varieties; the one approach, which is fairly standard, is through a system of semigroups as done for example in the first author's Ph.D. thesis. The other is by viewing a not-necessarily-normal toric variety as a torus orbit closure in some smooth ambient variety with a linear torus action (here, the meaning of ``linear'' can be made precise by embedding a smooth variety into the multi-Proj of a polynomial ring; see e.g., \cite{multi_proj}). The third, more abstract one, is to characterize toric varieties in terms of maps to them; i.e., functorial points of them (the normal case of this is due to Cox \cite{Cox_functor} and Kajiwara \cite{Kajiwara1998}, and the general case due to this paper as well as \cite{multi_proj}).

Now, here is the precise statement for the map $\varphi$.

%\begin{samepage}
\begin{thm}[Theorem \ref{thm:etale map to a special fiber}] 
Suppose we are given a toric degeneration $\mathfrak{X}$ from $X$ (a general fiber) to $W$ (the special fiber) over $\mathbb{C}$ that is embedded in a complex smooth variety as a proper family.

Then there exists a continuous surjective map (not a morphism of varieties)
$$\varphi : X \to W$$
such that, for each torus orbit $O \subset W$, $\varphi^{-1}(O)$ has the structure of a reducible variety so that $\varphi : \varphi^{-1}(O) \to O$ is a finite \'etale covering.

(Note we are not claiming that $\varphi^{-1}(O)$ is a subscheme of $\mathfrak{X}$.)
\end{thm}

\subsection{Symplectic consideration}\label{sec:symplectic considerations} One of the motivations for the construction of $\varphi$ above was an application to sympletic topology, as was originally done by Harada and Kaveh.

%This paper includes an application to symplectic geometry. This is in part because our construction of a map in the previous section was inspired by the previous work of Harada and Kaveh. They have used their map to a toric variety to construct an integrable system. Using our map, their construction can be repeated but with an important improvement.

First we recall the fact that a projective toric variety has a natural Hamiltonian group action relative to a choice of an embedding into a projective space. In particular, it comes with the moment map $\mu_W$, a stratum-preserving surjective map from a toric variety to a polytope, where the polytope is stratified by open faces. If $W$ is a special fiber of a toric degeneration of $X$ constructed by the Gr\"obner theory, then the image of the moment map is exactly the Newton–Okounkov body $\Delta$ of $X$. We can then simply compose the two maps:
$$\mu_X : X \overset{\varphi}\to W \overset{\mu_W}\to \Delta.$$
This map is the moment map for $X$ in the following sense:

%The construction of an integral system using a toric degeneration was a result of Harada and Kaveh. The above theorem extends their result to:

\begin{thm}[Harada--Kaveh type Theorem \ref{thm:improved Harada-Kaveh}] 
Let $\mathfrak{X}\subset M$ be a toric degeneration of a variety $X$ embedded in some complex smooth variety $M$ in such a way that $\mathfrak{X} \to \mathbb{A}^1$ is the restriction of a proper morphism $M \to \mathbb{A}^1$. 
Then there exists a continuous surjective map
$$\varphi : X \to W$$
that is stratum-wise étale in the sense in Theorem~\ref{thm:etale map to a special fiber} above.
Moreover, given the moment map $\mu_W$ of $W$ we have that
$$\mu_X = (f_1, \dots, f_r) : X \overset{\varphi}\longrightarrow W \overset{\mu_W}\longrightarrow \mathbb{R}^r$$
is an integrable system in the stratified sense (Definition \ref{integrable system in stratified}), i.e. linear combinations of $f_1, \dots, f_r$ form an integrable system on each stratum of $X$.
\end{thm}

The map $\varphi$ in Theorem~\ref{thm:etale map to a special fiber} may not preserve a symplectic structure, so we modify it by Weinstein's theorem (Theorem \ref{Weinstein}), which says the symplectic structure near $W$ is unique up to diffeomorphism, and the result is the above map. A key extension here over \cite{HK15} is that the integrable system is constructed also on the boundary, in addition to the fact that $\varphi$ is defined explicitly not as a continuous extension.

%As there is no standard reference on symplectic structures on singular spaces, Appendix gives a survey on approaches to symplectic structures in the singular situation, including ours.

%in Theorem~\ref{thm:embedding toric subalgebra} we construct such a (global) map explicitly in the case that the toric degenerations is \emph{equivariant}.
%that is the torus action on the generic fiber is extended equivariantly to the family (see Remark~\ref{rmk:equivariant TD}).
%As a consequence the toric degeneration may be realized as a ).

\subsection{Possible alternatives }
For the interested readers, we want to record a few speculations.

\begin{remark} In reviewing moment maps, we have included the convexity result over an algebraically closed field $k$ of characteristic zero (well known over $\mathbb{C}$). According to the experts, the Thom--Mather theory holds over a real closed field as well and if we assume that, the results of \S \ref{sec:collapsing a family} are also true over a real closed field. Hence, we could conclude that in characteristic zero, given a projective toric degeneration, the Newton--Okounkov body over a real closed field $k_{\mathbb{R}}$ is the image of a map to the special fiber followed by the moment map (this generalizes the known case $k_{\mathbb{R}} = \mathbb{R}$). This is moderately useful since for example, in representation theory, one might like to work with an algebraically closed field of characteristic zero not just $\mathbb{C}$.

Sympletic structures in real algebraic geometry seem not to have been studied much. However, the only part in this paper where we need a non-algebraic result is Weinstein's theorem; the proof here uses ordinary differential equations. An algebraic proof of it does not seem to be known. We speculate that one way to give an algebraic proof might be to adapt (to be precise, reverse) a quantized contact transformation used in the theory of microdifferential operators (cf. \cite[Ch. I., \S 5.]{schapira}), or at least ideas from that circle.
\end{remark}

\begin{remark} After the completing the paper, we learned about a paper of Drinfeld \cite{drinfeld} in which he generalizes the Białynicki-Birula decomposition to a singular space. The construction in that paper may give an alternative approach to this paper.
\end{remark}

\subsection{Further topics}\label{sec:further}

A special case of a map to a toric variety is a torus-equivariant vector bundle over a toric variety, the topic studied by, for example, C. Manon and K. Kaveh \cite{KM-toricbundles}.  Hence, the natural problem is to classify maps to toric varieties; in fact, the subsection at the end of the paper already gives a hint.
For the importance of the case of topological coverings, see a discussion around Theorem \ref{topological coverings of a complex variety}. (Very roughly speaking, our paper suggests there could be some Galois-type classification but in the stratified sense.)

%To take this further it is necessary to prove the local invariant cycle theorem for non-semistable degenerations.

A degeneration considered in literature is typically \emph{semistable} (i.e., a special fiber is a divisor with normal crossing) or at least reduced to such by semistable reduction. One such instance is the local invariant cycle theorem. Clemens proved the local invariant cycle theorem for semistable degenerations (that is, the special fiber is a divisor with normal crossings) in \cite{Clemens} by constructing a map from a general fiber to the special fiber.

The result of \S \ref{sec:collapsing a family} is notable in that it applies to a toric degeneration, which is \emph{not} semistable. It seems reasonable that Clemens' result extends to a non-semistable case by using the map in \S \ref{sec:collapsing a family}. 
Even more generally, nowadays the local invariant cycle theorem can be obtained as an immediate consequence of the {\it BBG decomposition} --the decomposition theorem of Beilinson, Bernstein and Deligne. 
As the main application of the construction in \S6 of \cite{MacPhersonGoresky} is to the specialization of perverse sheaves, it is probably possible to do a version of the BBG decomposition in the setup of \S\ref{sec:collapsing a family}.

\subsection{Acknowledgements}
We are grateful to Eduardo González for a detailed reading of an early version of the paper and his valuable feedback.
L.B. acknowledges support of the PAPIIT projects IA100122 (resp. IA100724) dgapa UNAM 2022 (resp. 2024) and CONAHCyT project CF-2023-G-106.

\section{Rees algebras in several parameters}\label{sec:Rees algebra}

As noted in the introduction, it is the most common to consider a flat degeneration over the affine line $\mathbb{A}^1$. However, in order to construct a torus-invariant stratification, we need to consider a flat degeneration over an affine space $\mathbb{A}^r$ or, algebraically, a Rees algebra in several parameters. 
Thus, in this section, we give some discussion on a multi-graded, generalized and extended Rees algebra, including the construction of a stratification (Construction \ref{locally Grobner toric degeneration}).

\

First, we note some facts on multi-graded rings and modules. If $M$ is a $\mathbb{Z}^r$-graded module over a $\mathbb{Z}^r$-graded ring $A$, then we write
$$\operatorname{sp}(M) = \{ a \in \mathbb{Z}^r \mid M_a \ne 0 \}.$$
The name sp stands for spectrum. We order $\mathbb{Z}^r$ by a lexicographical ordering and then we shall write $\mathbb{Z}^r_{\ge a} = \{ b \in \mathbb{Z}^r \mid b \ge a \}$ and similarly $\mathbb{Z}^r_{> a}$, $\mathbb{Z}^r_{\le a}$, etc.

\begin{lemma}[graded Nakayama]\label{graded Nakayama} Let $A$ be a $\mathbb{Z}^r$-graded ring and $I \subset A$ a graded ideal such that $\operatorname{sp}(I) \subset \mathbb{Z}^r_{> 0}$. If $M$ is a $\mathbb{Z}^r$-graded $A$-module such that $IM = M$ and $\operatorname{sp}(M)$ has a least element, then $M = 0$.
\end{lemma}
\begin{proof} Assume $M \ne 0$. We have:
$$(I M)_a = \bigoplus_{p > 0} I_p M_{a - p}.$$
Since $IM = M$, we get $\inf(\operatorname{sp}(M)) = \inf(\operatorname{sp}(IM)) > \inf(\operatorname{sp}(M))$, a contradiction.
\end{proof}

The next lemma contains the Noether normalization lemma as a consequence of the above graded Nakayama lemma.

\begin{lemma}\label{graded Noether normalization} Let $A$ be a $\mathbb{Z}^r$-graded ring such that $\operatorname{sp}(A)$ is a well-ordered subset of $\mathbb{Z}^r_{\ge 0}$. Assume $A_0$ is a Noetherian ring. Then the following are equivalent.
\begin{enumerate}
    \item $A$ is a finitely generated $A_0$-algebra.
    \item The ideal $A_+ = \bigoplus_{a > 0} A_a$ is finitely generated.
    \item For some homogeneous elements $y_1, \dots, y_d$ in $A$ such that $\sqrt{A_+} = \sqrt{(y_1, \dots, y_d) A}$, we have that $A$ is finite over $A_0[y_1, \dots, y_d]$.
\end{enumerate}
\end{lemma}
\begin{proof} A proof is well known when $r=1$ and we shall repeat the same proof.

(iii) $\Rightarrow$ (i) Clear. (i) $\Rightarrow$ (ii): $A_+$ is finitely generated if (i) is true, since $A$ is then a Noetherian ring.

(ii) $\Rightarrow$ (iii) Lemma \ref{graded Nakayama} gives the criterion: $M$ is finite over $A$ if and only if $M/IM$ is finite over $A/I$ for some graded ideal $I$ with $\operatorname{sp}(I) \subset \mathbb{Z}^r_{ > 0}$ (\emph{proof}: let $N \subset M$ be a finitely generated submodule which coincides with $M$ modulo $IM$. Then $N + IM = M$ and so $I(M/N) = M/N$, which implies $M/N = 0$.) 

Thus, to show $A$ is finite over $B := A_0[y_1, \dots, y_d]$, it is enough to show $A/IA$ is finite over $B/I$ for some graded ideal $I \subset B$ with $I \subset B_+$. We take $I = B_+$. First, note that $A/IA$ is integral over $B/I = A_0$. Indeed, $(A/IA)_0 = A_0$ and if $f$ is in $A$ with positive degree, then $f^n$ is in $IA$ for some $n > 0$ since $A_+ \subset \sqrt{A_+} = \sqrt{IA}$.

Second, we shall note that $A$ is of finite type (i.e., a finitely generated algebra) over $A_0$. For that, let $g_1, \dots, g_m$ be homogeneous generators of $A_+$ of multi-degrees $d_i$'s. Consider $E = \{ \alpha \subset \operatorname{sp}(A) \mid \bigoplus_{a \in \alpha} A_a \subset A_0[g_1, \cdots, g_m] \}$. Then, by Zorn's lemma, this set $E$ contains a maximal element $\alpha$. Suppose $\alpha \subsetneq \operatorname{sp}(A)$. Then, by well-ordering, the set $\operatorname{sp}(A) - \alpha$ contains a least element $a$. But
$$A_a = \sum_i A_{a-d_i} g_i$$
and $a - d_i$ is in $\alpha$ unless $A_{a - d_i} = 0$. This implies $A_a \subset A_0[g_1, \cdots, g_m]$, a contradiction. Hence, $\alpha = \operatorname{sp}(A)$.

Finally, $A/IA$ is integral and of finite type over $B/I$; thus, is finite over the latter (in general, finite = integral + of finite type).
\end{proof}

\begin{corollary} The above ring $A$ is a Noetherian ring if and only if it is a finitely generated $A_0$-algebra.
\end{corollary}

Now, let $R$ be a ring. We consider a decreasing filtration $\mathfrak{i}_* = \{ \mathfrak{i}_a \subset R \, \textrm{ subgroup} \mid a \in \mathbb{Z}^r \}$ of $R$ over $\mathbb{Z}^r$, where ``decreasing'' means $\mathfrak{i}_b \subset \mathfrak{i}_a$ for $b > a$. We say such a decreasing filtration $\mathfrak{i}_*$ is an \emph{ideal filtration} if 
\begin{enumerate}[label=(\arabic*)]
    \item $\mathfrak{i}_a = R$ for $a \le 0$,
    \item $\mathfrak{i}_a \mathfrak{i}_b \subset \mathfrak{i}_{a + b}$ for each $a, b$. 
\end{enumerate}
Note each $\mathfrak{i}_a$ is then an ideal of $R$ since $\mathfrak{i}_a R = \mathfrak{i}_a \mathfrak{i}_0 \subset \mathfrak{i}_a$. For notational simplicity, we often also write $\mathfrak{i}$ for $\mathfrak{i}_*$.

The notion has been around in commutative algebra for some time but the term ``ideal filtration'' is from the first author's Ph.D. thesis, while some authors simply use the term ``filtration''. A basic example is given by powers of an ideal; i.e., $\mathfrak{i}_n = I^n$. Another example is given by a valuation; namely, when $R$ is an integral domain and $\nu$ is a valuation on the field of fractions of $R$ that is $\ge 0$ on $R - 0$, we let $\mathfrak{i}_a = \{ 0 \} \cup \{ f \in R - 0 \mid \nu(f) \ge a \}.$

Given an ideal filtration $\mathfrak{a}$, a basic question is when the (extended) Rees algebra $\bigoplus_a {\mathfrak{a}}_a$ is finitely generated as an $R$-algebra. One answer is given by Proposition \ref{finite generation of a Rees algebra} below. The proposition involves the notion of an associated $\mathbb{Z}^r$-graded ring.

\begin{definition}[associated ring] For an ideal filtration $\mathfrak{a}_*$, we write
$$\operatorname{gr}_{\mathfrak{a}} R = \bigoplus_a \mathfrak{a}_a/\mathfrak{a}_{> a}$$
where $\mathfrak{a}_{> a} := \sum_{b > a} \mathfrak{a}_b = \bigcup_{b > a} \mathfrak{a}_b$. Note $\operatorname{sp}(\operatorname{gr}_{\mathfrak{a}} R) \subset \mathbb{Z}^r_{\ge 0}$.

For each $f$ in $R - \cap_a \mathfrak{a}_a$, we write $\sigma_{\mathfrak{a}}(f)$ or just $\sigma(f)$ for the image of $f$ in $\mathfrak{a}_a/\mathfrak{a}_{> a}$ when $f$ is in $\mathfrak{a}_a - \mathfrak{a}_{> a}$. It is called the \emph{initial form} of $f$; also called the \emph{principal symbol} of $f$ especially in the non-commutative context. Note, in general, $\sigma$ is not a group homomorphism for multiplication nor addition. It is just a set-theoretic map. We call the map $\sigma$ the \emph{symbol map}.
\end{definition}

The next proposition in a weaker form is often stated in the context of complete rings; cf. \cite[Proposition 7.12.]{Eisenbud}. (For a non-commutative example, there is \cite[Proposition 1.1.8.]{schapira}.)

\begin{prop}\label{finite generation of a Rees algebra} Let $R$ be a Noetherian ring and $\gamma \subset R$ a subset that generates $R$ as a ring or as a $k$-algebra if $R$ contains a field $k$.

Let $\mathfrak{a}$ be an ideal filtration on $R$ such that each $\mathfrak{a}_a$ is generated by a subset of $\gamma$ and $\mathfrak{a}_{> 0} \ne R$.

Assume for each increasing sequence $a_1 < a_2 < \cdots$ and each ideal $I$ generated by a subset of $\{ xy \mid x, y \in \gamma \}$, we have
$$I = \bigcap_i (I + \mathfrak{a}_{a_i}).$$

Now, if $\operatorname{gr}_{\mathfrak{a}} R$ is generated as an $R/\mathfrak{a}_{> 0}$-algebra by the image of a finite subset $\gamma_0$ of $\gamma$ under the symbol map $\sigma$ and if the semigroup $\mathfrak{S}$ generated by $\operatorname{sp}(\operatorname{gr}_{\mathfrak a} R)$ is well-ordered as a set, then the two rings
$$\bigoplus_{a \in \mathfrak S} \mathfrak{a}_a$$
and 
$$\mathfrak{R}(\mathfrak a) := \bigoplus_{a \in \mathbb{Z}\mathfrak{S}} u^{-a} \mathfrak{a}_a, \, u^b = u_1^{b_1} \dots u_r^{b_r},$$
are generated as $R$-algebras by
$$\{ g \, \textrm{ at } \deg(\sigma(g)) \mid g \in \gamma_0 \} \textrm{ and } \{ u^{-\deg(\sigma(g))} g \mid g \in \gamma_0 \}$$
respectively. Here, $\mathbb{Z}\mathfrak{S}$ is the group generated by $\mathfrak S$, $u_i$'s indeterminates and $\mathfrak{R}(\mathfrak a) \subset R[u_1^{\pm}, \dots, u_r^{\pm}]$ is a graded subring.
\end{prop}
\begin{proof} We have $$\bigoplus_{a \in \mathfrak S} \mathfrak{a}_a \simeq \bigoplus_{a \in \mathfrak S} u^{-a} \mathfrak{a}_a$$ by mapping each $f$ in the $a$-th component to $u^{-a} f$. We shall show the latter $A := \bigoplus_{a \in \mathfrak S} u^{-a} \mathfrak{a}_a$ is finitely generated as an $A_0 = R$-algebra with generators in $\gamma$. In light of Lemma \ref{graded Noether normalization}, we shall show that the ideal $A_+$ is finitely generated.

Let $g_1, \dots, g_m$ be elements of $\gamma$ such that $\sigma(g_i)$ are defined and generate the ideal $(\operatorname{gr}_{\mathfrak{a}} R)_+$. Let $d_i = \deg(\sigma(g_i))$. Then, for each $a > 0$ in $\mathbb{Z}^r$,
$$\mathfrak{a}_a = \sum \mathfrak{a}_{a - d_i} g_i + \mathfrak{a}_{> a}.$$
Since $R$ is Noetherian, $\mathfrak{a}_{> a} = \mathfrak{a}_{a_2}$ for some $a_2 > a$. Then, by the above applied to $\mathfrak{a}_{a_2}$, we get:
$$\mathfrak{a}_a = \sum \mathfrak{a}_{a - d_i} g_i + \mathfrak{a}_{> a_2}$$
since $\mathfrak{a}_{a_2 - d_i} \subset \mathfrak{a}_{a - d_i}.$ Let $I = \sum_i \mathfrak{a}_{a - d_i} g_i$. Thus, inductively, we get:
$$\mathfrak{a}_a = \bigcap_j (I + \mathfrak{a}_{a_j}) = I$$
where the equality on the right holds by assumption. That is,
$$u^{-a} \mathfrak{a}_a = \sum u^{-a + d_i} \mathfrak{a}_{a - d_i} u^{-d_i} g_i$$
and so $A_+$ is finitely generated.

Finally, $\mathfrak{R}(\mathfrak a)$ is an algebra over $R[u^{d_1}, \dots, u^{d_m}]$ that is finitely generated by $u^{-d_i} g_i$'s.

%and is in fact a finitely generated one since $A = \bigoplus_{a \ge 0} u^{-a} \mathfrak{a}_a$ is finitely generated as an $R$-algebra. Hence, $\mathfrak{R}(\mathfrak a)$ is a finitely generated algebra over $R$.
\end{proof}

\begin{remark} If the conclusion of the proposition holds, that is, $\mathfrak{R}(\mathfrak a)$ is finitely generated as an $R$-algebra, then we also have that $R$ is finitely generated as a $\mathbb{Z}$-algebra or a $k$-algebra if a field $k \subset R$.
\end{remark}

A trivial case in which the assumption $I = \cap_i (I + \mathfrak{a}_{a_i})$ above holds is when $\mathfrak{a}$ is nilpotent; i.e., $\mathfrak{a}_{a_i} = 0$ for some $i$. The next corollary says the degree-wise such nilpotency is enough.

\begin{corollary}\label{graded corollary} Assume $R$ is $\mathbb{N}_0$-graded and assume $\mathfrak{a}_a$ is graded as an ideal for each $a$. Also, assume for each $n$ and for each sequence $\mathfrak{a}_i$ in the proposition, there exists some $i$ such that $\mathfrak{a}_{a_i, n} = 0$.

Now, if $\operatorname{gr}_{\mathfrak{a}} R$ is finitely generated as an $R/\mathfrak{a}_{> 0}$-algebra, then $\mathfrak{R}(\mathfrak{a})$ is finitely generated as an $R$-algebra.
\end{corollary}
\begin{proof} Take $\gamma$ to be the set of all homogeneous elements in $R$. We shall check the condition
$$I = \bigcap_i (I + \mathfrak{a}_{a_i}).$$
It is enough to check this degree-wise and then the equality is clear by assumption.
\end{proof}

\begin{corollary}\label{valuation gr R} Let $R$ be a $\mathbb{N}_0$-graded Noetherian integral domain with $R_0 = k$ a field and of Krull dimension $r$. Then let
$$\nu : R - 0 \to \mathbb{Z}^{1+r}$$
be a valuation such that $\nu \ge 0$ and $\nu(k - 0) = 0$. Also, assume we can write $\nu = (\deg, \nu')$ for some graded valuation $\nu'$. Let $\mathfrak v$ be the ideal filtration associated to $\nu$; i.e., $\mathfrak{v}_a = \nu^{-1}(\mathbb{Z}^r_{\ge a}) \cup 0$.

If $\nu(R - 0)$ is finitely generated as a semigroup and the rational rank of the semigroup $\nu(R - 0)$ is $r + 1$, then $\mathfrak{R}(\mathfrak{v})$ is finitely generated as an $R$-algebra.
\end{corollary}
\begin{proof} Note $\gr_{\mathfrak v}R$ is the semigroup algebra of $\nu(R - 0)$ and so the assumptions in the preceding Corollary \ref{graded corollary} all hold.
\end{proof}

\begin{remark}[construction of a toric degeneration] It is still an open question when we can find a valuation $\nu$ in Corollary \ref{valuation gr R} (constructing $\nu$ is relatively easy and was done for example in \cite{Okounkov} but constructing $\nu$ such that $\nu(R - 0)$ is finitely generated is far more difficult). Aside from many known special cases, known general constructions include that in \cite{Anderson} and that in the first author's Ph.D. thesis. See also for example \cite{Rational_curve_example} for a worked-out example illustrating the difficulty.
\end{remark}

In the setup of Corollary \ref{valuation gr R}, let $S = k[u_1, \dots, u_r]$ be the polynomial ring and $\mathfrak{R} = \mathfrak{R}(\mathfrak v)$. Geometrically, we have
$$p : \Spec(\mathfrak{R}) \to \mathbb{A}^r = \Spec(S),$$
which is of finite type. Alternatively, we can take Proj with respect to the grading inherited from $R$
$$p : \Proj(\mathfrak{R}) \to \mathbb{A}^r,$$
which again is of finite type.

Also, from such $p$, we can get a degeneration over $\mathbb{A}^1$. Namely, choose some $\mathbb{A}^1 \hookrightarrow \mathbb{A}^r$ such that $\mathbb{A}^1$ and $\mathbb{A}^r$ share the same origin. Then the fiber $p|_{\mathbb{A}^1}^{-1}(0)$ is the same as the fiber $p^{-1}(0)$.

\begin{comment}
The flatness of the above $p$ is a bit tricky in the several-parameter case; in the one-parameter case, the question reduces to torsion-free-ness, while in the several-parameter case, there is no similar reduction. One way to see it is to use the below when it applies; for example, when $R$ is the homogeneous coordinate ring of a projective variety.

\begin{prop}[cf. Exercise 20.14. in \cite{Eisenbud}] Let $A$ be a $\mathbb{Z}^r$-graded Notherian integral domain. Let $M$ be a graded module over $A$ such that each degree piece $M_a$ is a finite $A$-module.

Then $M$ is flat if and only if the multi-Hilbert function of $M$ is constant over $\Spec(A)$; i.e., $\dim_{k(\mathfrak{p})} (M_a \otimes k(\mathfrak{p}))$ is independent of prime ideals $\mathfrak{p}$ of $A$ for each $a$.
\end{prop}
\begin{proof} We recall (\cite[Exercise 20.13.]{Eisenbud}) that, since $A$ is a Noetherian integral domain, a finite $A$-module is flat if and only if the function $\Spec(A) \to \mathbb{Z}, \, \mathfrak{p} \mapsto \dim_{k(\mathfrak{p})} (M \otimes k(\mathfrak{p}))$ is constant. Thus,
$M$ is flat $\Leftrightarrow$ each direct summand of $M$ is flat $\Leftrightarrow$ the multi-Hilbert function of $M$ is constant over $\Spec(A)$.
\end{proof}
\end{comment}

The degeneration over $\mathbb{A}^r$ has an advantage to one over $\mathbb{A}^1$ in that the former has the natural $\mathbb{G}_m^r$-action. This is important in the below:

\begin{construction}[of a torus-invariant stratification]\label{locally Grobner toric degeneration} Given a possibly reducible variety $X$, there is an inductive procedure to construct a stratification. Namely, first we write $X$ as a disjoint union of $X_0 = X_{reg}$ and the complement $Z = X - X_0$. Then we write $Z$ as a union of $Z_0$ and the complement and so on. In this procedure, we can also impose some additional conditions so that $Z_0$ is open dense subset of $Z_{reg}$ (i.e., we remove more than just singular points). For example, one can impose Whitney's condition (b) discussed in the next section and then the resulting stratification is a Whitney stratification; this is precisely Mather's canonical stratification in \cite[Theorem 4.9.]{Mather2}.

This procedure can also be applied to the morphism $p : \mathfrak{X} \to \mathbb{A}^r$ discussed above as well. That is, we construct a stratification $\Sigma$ that refines the stratification $\{ \mathfrak{X} - p^{-1}(0), p^{-1}(0) \}$. Precisely, we construct a sequence $\Sigma_d \subset \Sigma_{d-1} \subset \cdots \subset \Sigma_0$ of pairwise disjoint sets, $d = \dim \mathfrak{X}$, such that
\begin{enumerate}[label=(\arabic*)]
    \item each element of $\Sigma_i$, called a \emph{stratum}, is a locally closed smooth subvariety of either $\mathfrak{X} - p^{-1}(0)$ or $p^{-1}(0)$,
    \item each stratum in $\Sigma_i$ has dimension $\ge i$,
    \item For $|\Sigma_i| := \bigcup_{A \in \Sigma_i} A$, we have $|\Sigma_0| = \mathfrak{X}$.
\end{enumerate}
Then $\Sigma := \Sigma_0$ is a stratification we seek. The construction is by descending induction on the subscript $i$. Thus, first let $\Sigma_d = \{ \mathfrak{X}_{reg} \cap U \}, U = \mathfrak{X} - p^{-1}(0)$ and then suppose we have constructed $\Sigma_i$. Let $Y$ be the largest smooth open subset (namely, the regular locus) in $(\mathfrak{X} - |\Sigma_i|) \cap U$. Similarly, let $Y'$ be the largest smooth open subscheme in $(\mathfrak{X} - |\Sigma_i|) \cap p^{-1}(0)$. Then Let $\Sigma_{i-1}$ be the union of $\Sigma_i$ and all the irreducible components of $Y, Y'$ having dimension $\ge i - 1$.

The above might look needlessly formal. But the point is that because the construction is procedural, we can impose additional conditions in each step; in practice, we would require the above $Y$ to be more than just smooth.

Also, we note that, because the procedure is torus-invariant, each stratum in $\Sigma$ is necessarily $\mathbb{G}_m^r$-invariant.

\end{construction}

\begin{remark} An attentive reader would have noticed that the above construction never uses projectivity or affineness. So, it is actually the most natural to state it for a toric degeneration that is \emph{locally Gr\"obner} in some sense. However, such a local Gr\"obner theory is an open problem as already remarked in the introduction in \cite{algebraic_tex}. Note a naive such notion in the Zariski topology does not seem to work well since the assumption on $I$ in Proposition \ref{finite generation of a Rees algebra} generally fails Zariski-locally. (Of course, one might simply consider some abstract $\mathfrak{X} \to \mathbb{A}^r$ with a torus action, but the torus action then may not be linear and we do not have a theory of toric varieties with non-linear torus action.)
%n fact, the remaing sections of this paper is 
\end{remark}

\begin{comment}

In the case of a toric degeneration, the above gives

\begin{prop} Let $p : \mathfrak{X} \to \mathbb{A}^1$ be a toric degeneration. Assume, Zariski-locally, $\mathfrak{X}$ is obtained from $\mathfrak{Y} \to \mathbb{A}^r$.

Let $(*)$ be a property for the pair of strata $(A, B)$ that is isomorphism-invariant.

Then there exists a stratification $\Sigma$ on $\mathfrak{X}$ with 
\end{prop}
\begin{proof}
\end{proof}
\end{comment}

\section{Collapsing map and its restriction}\label{sec:collapsing a family}

The purpose of this section is to construct the following: given a smooth complex variety $M$ and a Zariski-closed subset $W \subset M$, in the classical topology, we construct a not-necessarily-open neighborhood $W(\delta)$ of $W$ and a continuous map
$$\pi : W(\delta) \to W$$
that is a retract of $i : W \hookrightarrow W(\delta)$; i.e., $\pi \circ i$ is the identity map. Moreover, it has the properties, making it a substitute for a tubular neighborhood of $W$: for the interior $U$ of $W(\delta)$, 
\begin{enumerate}[label=(\arabic*)]
\item $\pi : U - W \to M$ is smooth and $\pi : U \to M$ is locally Lipschitz continuous.
\item $\pi|_U$ is a strong deformation retract of $i$ (i.e., there is a homotopy from it to the identity such that the homotopy is the identity on $W$).
\item A continuous version of the relative Poincar\'e lemma holds for $\pi|_U$ (see Remark \ref{continous forms}).
\end{enumerate}

The existence of a similar map is due to \cite[\S7]{Goresky} and \cite[\S6]{MacPhersonGoresky}. Our construction is relatively simpler and also we get a strong deformation retract instead of a deformation retract. The strong version is needed for our application on integrable systems in the next section. In a bit more details, their construction is based on the Thom--Mather stratification theory, which constructs a system of tubular neighborhoods of the strata on $W$. Roughly, composing the projections corresponding to the tubular neighborhoods gives the claimed $\pi$. The required properties then hold since they hold for the tubular neighborhood projections.
%In the usual formulation of the theory, the construction is topological, which is not enough. Here, we thus actually use a semialgebraic version of the theory known to the experts; this is where the algebraic structures on $M, W$ are used.

In the case of a degeneration $\mathfrak{X} \subset M$, the above $\pi$ leads to the following key result. Namely, we construct $\pi$ so that, in addition to having the above properties, $\pi$ is stratum-wise locally trivial on a general fiber of $\mathfrak{X}$.

\begin{theorem}\label{thm:etale map to a special fiber} Let $p : \mathfrak{X} \to \mathbb{A}^1$ be a surjective proper morphism from a complex algebraic variety.
Assume that $\mathfrak{X}$ is a closed subvariety of a smooth complex algebraic variety $M$ and that $p$ is the restriction of a proper morphism $M \to \mathbb{A}^1$; e.g., $M = \mathbb{P}^n \times \mathbb{A}^1$ and $p$ the projection.

Then, for some $0 \ne t \in \mathbb{A}^1$ and some Whitney stratification on $W$ (see Remark \ref{canonical stratification}), there exists a surjective continuous map
$$\varphi : p^{-1}(t) \to W$$
such that for each stratum $A$ of $W$, the pre-image $\varphi^{-1}(A)$ has an induced reducible-variety structure so that $\varphi : \varphi^{-1}(A) \to A$ is a finite \'etale morphism.

Here, ``induced" means induced by $\varphi$ not by $p^{-1}(t)$. That is, we are *not* claiming $\varphi^{-1}(A)$ is a subscheme of $p^{-1}(t)$.
\end{theorem}

\begin{remark} The morphism $p$ in the theorem is not assumed to be flat but it is actually flat, since $\mathfrak{X}$ is a variety and a torsion-free module over a principal ideal domain is flat (\cite[Corollary 6.3.]{Eisenbud}). Also, if $p$ is trivial away from the special fiber $W = p^{-1}(0)$, then the statement holds for any $t \ne 0$ since it holds for some $t \ne 0$.
\end{remark}

%In the case of a toric degeneration, the stratification will be taken to be the orbit decomposition (so a stratum is an orbit); see Remark \ref{canonical stratification}. 

%\begin{remark}[on the literature]\label{collpasing literature} According to the discussion in \cite{mathoverflow}, it is a folklore result that a map $\varphi$ as in Theorem~\ref{thm:etale map to a special fiber} exists, perhaps with less properties (for example, apparently it is known to Deligne) at least when a general fiber is smooth. In fact, when the general fibers are smooth, the comments in \cite{mathoverflow} suggest that a retraction from a neighborhood of the special fiber to the special fiber can be constructed using a flow, even when a special fiber is non-smooth. 
%One suggestion there is to consider a flow on the resolution of singularities of the family $\mathfrak{X}$ and then compose the flow with the map from the resolution to $\mathfrak{X}$. 
%Meanwhile, Harada and Kaveh in \cite{HK15} give the same result using the uniform continuity to do a continuous extension. 
%When $p$ is semistable (meaning the special fiber is a divisor with normal crossings) and the general fibers are smooth, the existence of such a collapsing map is a result of Clemens in \cite{Clemens}.

%The key difference with the construction of MacPherson and Goresky and that of us is that they and we do not use a limit or continuity argument and that the map respects the local structure of a family; these are key to extending the results of \cite{HK15} outside the regular locus.
%\end{remark}

\begin{remark} After completing the first draft of this paper, we noticed the paper \cite{pflaum2019equivariant}, which also constructs a strong deformation retract of a neighborhood of $W$ in $\mathfrak{X}$. Their approach is to use a flow and thus is different from ours. We note a flow approach does not generalize to a real-algebraic geometry setup like here. See also \cite{mathoverflow} for a related discussion.
\end{remark}

\begin{remark}[continuous forms]\label{continous forms} By a \emph{continuous $k$-form}, we mean a continuous (as opposed to smooth) section of the $k$-th exterior power of the cotangent bundle. Throughout the paper, we usually have to work with continuous forms as opposed smooth forms. For example, this is essential to state and prove a continuous version of the relative Poincar\'e lemma mentioned earlier (Corollary \ref{cor:Poincare}).
\end{remark}

\begin{remark}[Whitney stratification]\label{canonical stratification} For the definition of a Whitney stratification, we refer to \cite[Part 1. Ch. 1.]{StratifiedMorse} as well as \cite{Mather}.
It is known that a subvariety of a smooth complex algebraic variety $M$ admits a Whitney stratification into smooth complex algebraic subvarieties of $M$ \cite[Part 1, Theorem 1.7.]{StratifiedMorse} or Construction \ref{locally Grobner toric degeneration}. Such a stratification is then essentially canonical; see the discussion following \cite[Theroem 4.9,]{Mather2}).

It also is known that, when an algebraic group is acting on a complex algebraic variety with finitely many orbits, the orbit decomposition of the variety is a Whitney stratification \cite[Theorem 7.6]{MacPherson}.
%Thus, in particular, in the case when $p$ is a toric degeneration the strata in Theorem \ref{thm:etale map to a special fiber} can be taken to be torus-orbits.

Since we are only interested in complex algebraic varieties in this paper, by Whitney stratification, we will always mean a stratification into complex smooth algebraic varieties or stratifications induced by them, unless said otherwise.
\end{remark}

For the Thom--Mather theory, we need some notations and terminology. Given two strata $A, B$, we write $A \le B$ if $A$ lies in the closure of $B$. Two strata are said to be {\it incomparable} if neither $A \le B$ nor $B \le A$ holds.
By a {\it tubular neighborhood} of a submanifold $A$ in $M$, we mean a smooth embedding $\psi_A : E \to M$ from a smooth vector bundle $E$ over $A$ to $M$ that is the identity on $A$ and the image of $\psi_A$ is an (open) neighborhood of $A$ in $M$.

Given such a $\psi_A$, we write $\mathfrak{T}_A$ for the image of $\psi_A$ and $\pi_{A}$ for the composition $\mathfrak{T}_A \overset{\psi_A^{-1}} \to E \overset{\xi}\to A$ where $\xi$ is the natural projection. 
Note that $\mathbb{R}$ acts on $E$ by scalar multiplication on each fiber. Hence, $\pi_{A}$ is a strong deformation retract of $A \subset \mathfrak{T}_A$ (i.e., there is a homotopy from $\pi_A$ to the identity on $\mathfrak{T}_A$ such that the homotopy is the identity on $A$). 
We also put an inner product on each fiber of $E$ such that the inner products vary smoothly across fibers and then write $\rho_{A} : \mathfrak{T}_{A} \to \mathbb{R}$ given as $\rho_{A}(x) = |\psi_A^{-1}(x)|^2$ where $| \cdot |$ denotes the norm on the fiber passing through $\psi_A^{-1}(x)$.

In his notes \cite{Mather} (specifically Proposition 7.1. in op. cit.), Mather constructs a system of tubular neighborhoods inductively so that each tubular neighborhood of a stratum is compatible with those of the strata in the boundary. Precisely, by induction on the dimension of the strata, he constructs the system $\{ \psi_A \mid A \}$ of tubular neighborhoods of $A$ such that, for strata $A < B$, we have
$$\pi_A \circ \pi_B = \pi_A, \rho_A \circ \pi_B = \rho_A$$
where the equations are defined.
Shrinking $\mathfrak{T}_A$'s, we shall assume $\mathfrak{T}_A, \mathfrak{T}_B$ are disjoint if $A, B$ are incomparable \cite[(A11) after Definition 8.2.]{Mather}.

%\begin{remark}[Conically smooth stratified spaces] Recently, the notion of a conically smooth stratified space has been introduced by Ayala, Francis and Tanaka \cite{ConicallySmooth}. The above system of tubular neighborhoods as well as the conical structures induced by it (Proposition \ref{conical structure} below) then may be thought of as the structure of a conically smooth stratified space; this is claimed to be proved by \cite{Nocera} (the proof of \cite[Theorem 3.4.]{Nocera} is somehow unclear but Proposition \ref{conical structure} and Remark \ref{cone structure} here can be used as a fix.)
%\end{remark}

We need the following variant of Thom's second isotopy Lemma \cite[Proposition 11.2.]{Mather}:

\begin{lemma}\label{second isotopy lemma} Fix a stratum $A$ of $W$. Let $S$ be a Whitney-stratified locally closed subset of $\mathfrak{T}_A$ and $f : S \to N$ the restriction of a smooth map $M \to N$, $N$ a smooth manifold. Assume that $f$ is a proper submersion on each stratum. Then for each point $y$ in $N$, there exists a neighborhood $U$ of $y$ in $N$ and a homeomorphism $H$ such that the diagram 
$$\xymatrix{f^{-1}(y) \times U \ar[rd]_{\pi_A \times \operatorname{id}} \ar[rr]_{\sim}^*[@]{H} & & f^{-1}(U) \ar[ld]^{(\pi_A, f)} \\
& A \times U}\\$$
commutes, $H$ preserves the strata and its restriction to each stratum is a diffeomorphism onto the image.

Also, if $N = \mathbb{R}^n$, then $U$ can be taken to be $\mathbb{R}^n$.
\end{lemma}

There are two proofs: a direct proof and one adapting the proof of \cite[Theorem 8.6.]{Mather2}. For later reference, we include both proofs.

\begin{proof}[Direct proof of Lemma~\ref{second isotopy lemma}]
Since the statement is local, it is enough to prove the case $N = \mathbb{R}^n$ with $U = \mathbb{R}^n$. For the notational clarity, we shall also assume $n = 1$ (the proof for the general case is the same).

Shrinking the ambient manifold $M$, we can assume $S$ is a closed subset. Then, by Thom's first isotopy lemma \cite[Proposition 11.1.]{Mather}, we have a trivialization $H$:
$$\xymatrix{f^{-1}(y) \times U \ar[rd]_{\operatorname{id}} \ar[rr]_{\sim}^*[@]{H} & & f^{-1}(U) \ar[ld]^{f} \\
& U.}\\$$
Thus, we only need to verify that the asserted diagram commutes with the above $H$; i.e., $\pi_A(H(x, t))$ is independent of $t$. In the proof of the first isotopy lemma, the trivialization $H$ is given by the flow $\varphi_t$ generated by a set-theoretic vector field $\eta$ smooth on each stratum; namely, $H(x, t) = \varphi_t(x)$. 
By a \emph{set-theoretic vector field}, we mean a not-necessarily-smooth-nor-continuous section of a tangent bundle.
Now, by construction, $\eta$ is tangent everywhere to the fibers of $\rho_A$. Thus,
$$\frac{d}{dt} (\pi_A \circ \varphi_t) = d \pi_A (\eta \circ \varphi_t) = 0.$$
That is to say, $\pi_A \circ \varphi_t = \pi_A \circ H(\cdot, t)$ is independent of $t$. This proves the claim.
\end{proof}

\begin{proof}[Mather's proof of Lemma~\ref{second isotopy lemma} adapted]
We follow \cite[Proof of Theorem 8.6.]{Mather2}. The projection $p : A \times \mathbb{R}^n \to \mathbb{R}^n$ is trivial in the sense that we have the commutative diagram
$$\xymatrix{p^{-1}(y) \times \mathbb{R}^n \ar[d]_{\operatorname{pr}} \ar[r]^-{p \times \operatorname{id}} & y \times \mathbb{R}^n \simeq \mathbb{R}^n \ar[d]^{\parallel}\\
A \times \mathbb{R}^n \ar[r]^-{p} & \mathbb{R}^n.}$$

The hypothesis on $f$ implies that $f$ satisfies Thom's condition $(a_f)$ defined in \cite[\S 11]{Mather}. Indeed, if $B$ is a stratum of $S$, then the tangent space to $f|_B^{-1}(0)$ at $x$ equals $\operatorname{ker}(d(f|_B)_x)$ since $f|_B$ is a submersion; thus, Thom's condition holds by Whitney's condition (A) (which is weaker than Whitney's condition (B)).

Then, by an argument preceding Proposition 2.9. of \cite{Bekka}, we see that $(\pi_A, f) : S \to A \times \mathbb{R}^n$ is a Thom mapping (see \cite[\S 11]{Mather}).
Then, as in the proof of \cite[Theorem 8.6.]{Mather2}, we find a homeomorphism $H$ such that the following commutes:
\[
\xymatrix{f^{-1}(y) \times \mathbb{R}^n \ar[r]^-{\,\,\, (\pi_A, f) \times \operatorname{id}} \ar[d]_{H}^{\wr} & p^{-1}(y) \times \mathbb{R}^n \ar[d]_{\operatorname{pr}} \ar[r]^-{p \times \operatorname{id}} & \mathbb{R}^n \ar[d]^{\parallel}\\
S \ar[r]^{(\pi_A, f)} & A \times \mathbb{R}^n \ar[r]^-{p} & \mathbb{R}^n.}
\]
This is exactly what is asserted.
\end{proof}

The next proposition is Theorem 8.3. of \cite{Mather2} and it establishes the \emph{conical} structure of a Whitney stratification (``conical'' is defined below). In fact, given this result, some authors define a stratified space as a space that has a stratification with the conical structure; see e.g., \cite[Ch I., \S 1.]{Borel}.

\begin{prop}[Conical structure]\label{conical structure} Given a stratum $A$ of $W$, consider
\[
l_A:\mathfrak{T}_A\to \mathbb R, \quad \text{given by} \quad x\mapsto \frac{\rho_A(x)^{1/2}}{(\epsilon\circ\pi_A)(x)},
\]
where $\epsilon$ is a positive smooth function on $A$. If $\epsilon$ is small enough, then
$(\pi_A, l_A) : l_A^{-1}((0, 1]) \to A \times (0, 1]$ is trivial over $\mathbb{R}$ as a family over $A$ in the sense that there is a homeomorphism $H_A$ over $A \times (0, 1]$:
    \[
    \xymatrix{l_A^{-1}(1) \times (0, 1] \ar[rd]_{\pi_A \times \operatorname{id}} \ar[rr]_{\sim}^*[@]{H_A} & & l_A^{-1}((0, 1]) \ar[ld]^{(\pi_A, l_A)} \\
    & A \times (0, 1]}
    \]
such that $H_A$ preserves the strata and the restriction of it to each stratum is a diffeomorphism onto the image. (Cf. Remark \ref{Hironaka}.)
\end{prop}
\begin{proof} By Lemma 7.3. and the remark preceding it in \cite{Mather}, we can choose $\epsilon > 0$ such that $(\pi_A, l_A)$ is proper and also, for each stratum $B$ on $W$, it is a submersion on $l_A^{-1}([0, 1]) \cap B$. Then $l_A^{-1}((0, 1])$ has a Whitney stratification induced from that on $W$ (see the proof \cite[Theorem 8.3.]{Mather2}); namely, the strata are $l_A^{-1}((0, 1]) \cap B$, $B$ strata of $W$ and $l_A^{-1}([0, 1]) - W$. The assertion then follows from Lemma \ref{second isotopy lemma} applied to $f = l_A : S = l_A^{-1}((0, 1]) \to \mathbb{R}$.
\end{proof}

The reason why the above gives a conical structure is as follows.

\begin{corollary}\label{mapping cylinder} $H_A$ identifies $l_A^{-1}([0, 1])$ as the mapping cylinder of $\pi_A : l_A^{-1}(1) \to A$; i.e.,
there is a push-out diagram:
$$\xymatrix{l_A^{-1}(1) \ar[r]^-{(\operatorname{id}, 0)} \ar[d]_{\pi_A} & l_A^{-1}(1) \times [0,1] \ar[d]\\
A \ar@{^{(}->}[r] & l_A^{-1}([0, 1])}$$
where the right vertical map is a unique extension of $H_A$ (which we shall also denote by $H_A$).
%The unique extension of $H_A$ will also be denoted by $H_A$.
\end{corollary}
\begin{proof} Let us extend $H_A$ by setting $H_A(x, 0) = \pi_A(x)$. We shall show that $H_A$ is then continuous.

Since the problem is near $A$, we may assume $M$ is a coordinate chart such that, by the uniqueness of a tubular projection, $\pi_A$ is the orthogonal projection onto $A$ and $\rho_A(x) = |x - \pi_A(x)|^2$. Then, for $x = H_A(x_0, s)$, since $\pi_A(x) = \pi_A(x_0)$, we have:
$$|H_A(x_0, s) - H_A(x_0, 0)| = |x - \pi_A(x)| = \rho_A(x)^{1/2} \le c \, l_A(x) = c s$$
for some constant $c$. Also, with $x = H_A(x_0, t), y = H_A(y_0, t)$, 
\begin{align*}
|H_A(x_0, t) - H_A(y_0, t)| &\le |x - \pi_A(x)| + |\pi_A(x_0) - \pi_A(y_0)| + |\pi_A(y) - y| \\
&\le 2ct + |\pi_A(x_0) - \pi_A(y_0)|.
\end{align*}
The continuity of $H_A$ now follows easily.
\end{proof}

In fact, we can do further trivialization

\begin{remark}\label{cone structure} By the submersion theorem applied to $f = \pi_A|_{l_A^{-1}(1)}$, for each point $a \in A$, we can find a neighborhood $U$ of $a$ in $M$ and a diffeomorphism over $A$:
$$(U \cap A) \times L \simeq l_A^{-1}(1) \cap \pi_A^{-1}(U \cap A)$$
    where $L = \pi_A^{-1}(a) \cap l_A^{-1}(1)$, called the \emph{link} of $a$. For simplicity, replacing $A$ by $U \cap A$, assume $A \subset U$. Then the above means $A \times L \simeq l_A^{-1}(1)$, which, together with $H_A : l_A^{-1}((0, 1]) \simeq l_A^{-1}(1) \times (0, 1]$, gives the strata-preserving $A$-homeomorphism:
$$l_A^{-1}([0, 1]) \simeq A \times cL, \, x \mapsto (\pi_A(x), l_A(x) u(x))$$
where $u : l_A^{-1}((0, 1]) \to l_A^{-1}(1) \to L$. In general, the above trivialization holds locally.
\end{remark}

\begin{remark}\label{Hironaka} As noted after \cite[Theorem 8.3.]{Mather2}, the trivialization $H_A$ in the above proposition cannot be taken to be smooth in general. However, since $W$ is a complex algebraic variety, we can (and do) assume the strata on $W$ are complex-algebraic. Now, it is known that Thom's first and second isotopy lemmas hold for Nash manifolds; i.e. semialgebraic real-analytic manifolds \cite{coste_nash,Escribano}. In particular, we can assume the trivialization $H_A$ as well as the trivialization in Remark \ref{cone structure} are Nash and so smooth.

Alternatively, it is sometimes mentioned (though we could not find a precise statement) in literature that, from Hironaka's works, we can find the trivializations like the above that are real-analytic. Finally, \cite{Verdier} also provides a trivialization that is more than a homeomorphism, which might also work for us.
\end{remark}

Using the above Proposition \ref{conical structure}, we shall define an extension of $\pi_A$ from some subset of the tubular neighborhood $\mathfrak{T}_A$. First, by Whitney's theorem (which says a closed subset of a manifold is the zero set of a smooth function), we choose a smooth real-valued function $\delta$ on $M$ such that $\delta^{-1}(0) = W$ and $0 \le \delta \le 1/4$ identically (explicitly, if $W = f^{-1}(0)$ for some smooth function $f$, take $\delta$ to be a multiple of $p \circ f$ for some smooth bounded non-negative function $p$ on $\mathbb{R}$ such that $p^{-1}(0) = 0$.)

For a stratum $A$ of $W$, we write $A(\delta) = \{ x \in \mathfrak{T}_A \mid l_A(x) \le \delta(x) \}$. Let $\lambda_A$ be a smooth function on $M - W = \{ \delta > 0 \}$ such that $\lambda_A = 1$ on $A(\delta)$, $\lambda_A < 1$ on $\mathfrak{T}_A - A(\delta)$ and $\Supp(\lambda_A) \subset A(2\delta)$.

Then define $\widetilde{\pi}_A : M \to M$ by:
\[
\widetilde{\pi}_A(x):=\left\{
\begin{matrix}
    H_A(x_0,(1-\lambda_A(x))s), & \text{if } x \in A(2\delta) \cap \{ \delta > 0 \} \text{ and } x=H_A(x_0,s)\\
    x, & \text{otherwise}.
\end{matrix}
\right.
\]

Note $\widetilde{\pi}_A$ is then an extension of $\pi_A|_{A(\delta)}$. Also,

\begin{lemma} $\widetilde{\pi}_A : M \to M$ is a continuous map and $\widetilde{\pi}_A : M - \overline{A} \to M$ is a smooth map. However, $\widetilde{\pi}_A : M \to M$ is not smooth (not even continuously differentiable.)
\end{lemma}

\begin{proof} Since $H_A$ is smooth away from $W$ by construction, the smoothness is clear away from $W$. Thus, we only need to show that $\widetilde{\pi}_A$ is smooth on some neighborhood of each point $x$ of $W$. 
Since $\widetilde{\pi}_A$ is smooth (in fact, it is the identity) outside the tubular neighborhood $\mathfrak{T}_A$ and $A$ is closed in $\mathfrak{T}_A$, we can assume $A$ is closed.

If $x$ is in $W - A$, then let $U \subset \mathfrak{T}_A$ be a neighborhood of $x$ that is disjoint from $A$. Then we have $U \subset \{ l_A \ge c \}$ for some constant $c > 0$. Since $\delta(x) = 0$, shrinking $U$, we can assume $U \subset \{ \delta \le c/2 \}$. But then $l_A / \delta \ge 2$ on $U$. Thus, $\lambda_A = 0$ on $U$ and so $\widetilde{\pi}_A$ is the identity on $U$, in particular smooth.

Next, we shall show $\widetilde{\pi}_A$ is continuous. For that, given the above, we only need to show $\widetilde{\pi}_A$ is continuous at each point $x$ of $A$. We write $x = H_A(x_0, 0)$. For any $y$ near $x$, if $y$ is in $\{ \delta > 0 \} \cap A(2\delta)$, then we write $y = H_A(y_0, s)$ and then $\widetilde{\pi}_A(y) = H_A(y_0, (1 - \lambda_A(y))s)$. Here, $|1 - \lambda_A(y)| \le 1$ and thus $(y_0, (1 - \lambda_A(y))s) \to (x_0, 0)$ as $(y_0, s) \to (x_0, 0)$. Otherwise, $\widetilde{\pi}_A(y) = y$. In any case, we have $\widetilde{\pi}_A(y) \to x = \widetilde{\pi}_A(x)$ as $y \to x$. (See Lemma \ref{Lipschitz} for a stronger result.)

Finally, we show that the derivative $\widetilde{\pi}_A'$, if it exists at all, is not continuous. For simplicity, via the trivialization $H_A$, view $\widetilde{\pi}_A$ as a map defined on $(0, 1] \times l_A^{-1}(1)$. The Jacobian matrix of it in the local coordinates $(y_1, \dots, y_n, s)$ is then
$$\begin{bmatrix}
1 & & 0 \\
\frac{\partial \lambda_A}{\partial y}s & & 1 - \lambda_A(y, s) - \frac{\partial \lambda_A}{\partial s} s
\end{bmatrix}
$$ where $\frac{\partial \lambda_A}{\partial y}s$ denotes the matrix consisting of partial derivatives. Here, $\lambda_A(y, s)$ does not go to $0$ as $s \to 0$. That is, the derivative $\widetilde{\pi}_A'$ represented by the Jacobian matrix is not continuous.
\end{proof}

\begin{remark} Here is a more conceptual reason why $\widetilde{\pi}_A$ is not smooth. If $\widetilde{\pi}_A$ \emph{were} smooth, a discussion later implies that there is a smooth retraction $U \to U$ of $W \hookrightarrow U$ for some neighborhood $U$ of $W$ in $M$. But, since $U$ is a manifold, that implies $W$ is a manifold (\cite[Ch 1. \S 2, Exercise 2.]{hirsch}), which is usually not the case.
\end{remark}

We still, however, have the following stronger continuity property, which is essential in \S \ref{sec:moment maps}.

\begin{lemma}\label{Lipschitz} With an appropriate choices of $\lambda_A$, we have that $\widetilde{\pi}_A : M \to M$ is locally Lipschitz continuous; i.e., each point of $M$ has a neighborhood $U$ in some coordinate chart together with a constant $c_U$ such that for all $x, y \in U$,
$$|\widetilde{\pi}_A(x) - \widetilde{\pi}_A(y)| \le c_U |x - y|$$
where $| \cdot |$ denotes the Euclidean length coming from the coordinate chart.
\end{lemma}

(In the paragraph just prior to \cite[Theorem 1.3.]{teschl}, ``locally Lipschitz continuous'' is defined so that the estimate holds on each compact set but it is known and is easy to see that that condition is equivalent to the condition given here.)

\begin{proof} It is enough to prove the assertion for a neighborhood $U$ of a point in $A$ such that $\overline{U}$ is compact and contained in some coordinate chart. Shrinking $U$, we shall assume $U \subset l_A^{-1}([0, 1]).$

Let $x = H_A(x_0, s), y = H_A(y_0, t)$ in $U \cap \{ \delta > 0 \}$. We have:
$$|\lambda_A(x) s - \lambda_A(y) t| \le s |\lambda_A(x) - \lambda_A(y)| + |s - t|.$$
Take $\lambda_A = q \circ (l_A/\delta)$ where $q$ is a smooth function with $q = 1$ on $[0, 1]$, $< 1$ on $(1, 2]$ and $= 0$ on $[2, \infty)$. By the mean value inequality,
\begin{align*}
s |\lambda_A(x) - \lambda_A(y)| &\le \sup|q'| \, s \, |s/\delta(x) - t/\delta(y)| \\
&\le \sup |q'| \left( \frac{s}{\delta(x)} |s - t| + \frac{st}{\delta(x)\delta(y)} |\delta(x) - \delta(y)| \right).
\end{align*}
Now, we can assume $\lambda_A(x), \lambda_A(y)$ are not both zero; without loss of generality, assume $x \in \operatorname{Supp(\lambda_A)} \subset A(2\delta)$. Then $s/\delta(x) \le 2$. If $y$ is also in $\Supp(\lambda_A)$, then the above is bounded by $\sup|q'| \, (2 |s - t| + 4 |\delta(x) - \delta(y)|).$

Assume $y \not\in \Supp(\lambda_A)$. Then let $c = t/2$ so that $t/c = 2$.  Since $\lambda_A(y) = 0 = q(t/c)$, the same analysis as above gives 
$$s |\lambda_A(x) - \lambda_A(y)| \le \sup|q'| \, s \, |s/\delta(x) - t/c|.$$
Here, since $s \delta(y) \le s c \le t \delta(x)$,
\begin{align*}
s \left| \frac{s}{\delta(x)} - \frac{t}{c} \right| &= \frac{s}{\delta(x) c} | s c - t \delta(x) | \le \frac{2}{c} | s \delta(y) - t \delta(x) | \le \frac{2\delta(y)}{c} | s - t | + \frac{2t}{c} | \delta(y) - \delta(x)|\\
&\le 2 | s - t | + 4 | \delta(y) - \delta(x)|.
\end{align*} Overall, we have:
$$\textstyle s |\lambda_A(x) - \lambda_A(y)| \le \sup|q'| \, (2 |s - t| + 4 \sup_{\overline{U}} |\delta'| \, |x - y|).$$

Now, since $1/({\epsilon \circ \pi}_A)$ is smooth and since $\rho_A^{1/2}$ is Lipschitz continuous as in the proof of Corollary \ref{mapping cylinder}, we have:
$$|s - t| = |l_A(x) - l_A(y)| \le c_{l_A}|x - y|$$
for some constant $c_{l_A}$. Hence, for some constant $c_{\lambda_A}$,
$$|\lambda_A(x) s - \lambda_A(y) t| \le c_{\lambda_A}|x - y|.$$

Next, define $\varphi_t : l_A^{-1}([0, 1]) \to l_A^{-1}([0, 1]), \, 0 \le t \le 1$ by $\varphi_t(H_A(x_0, s)) = H_A(x_0, ts)$. Fix a point $a_0 \in A$. Then Remark \ref{cone structure} gives $l_A^{-1}([0, 1]) \simeq A \times cL$ near $a_0$. 
Choosing coordinates on $L$, we have $L \subset \textrm{$l$-sphere} \subset \mathbb{R}^{l+1}$. Thus, again near $a_0$, we have $l_A^{-1}([0, 1]) \subset A \times \mathbb{R}^{l+1}$ where the inclusion preserves the strata and is smooth on each stratum. Then we see $\varphi_t$ can be given as $\varphi_t(a, u) = (a, tu)$ near $a_0$; in particular, $\varphi_t$ is Lipschitz continous near $a_0$ uniformly with respect to $t$.

Having the above estimates, it is now clear that $\widetilde{\pi}_A$ is locally Lipschitz continuous (note $\widetilde{\pi}_A(x) = \varphi_{(1 - \lambda_A(x))l_A(x)}(x)$.)
\end{proof}

By definition, the \emph{height} of a stratum $A$ is the supremum of the lengths $k$ of the chains $A = A_k > A_{k-1} > \cdots > A_0$, which are all finite since $A > B$ implies $B$ has strictly smaller dimension than $A$ \cite[Proposition 2.7.]{Mather}.

For each integer $i$, define $\widetilde{\pi}_i : M \to M$ by $\widetilde{\pi}_i|_{A(2\delta)} = \widetilde{\pi}_A|_{A(2\delta)}$ for the strata $A$ of height $i$ and $\widetilde{\pi}_i = $ the identity elsewhere (this is possible since our tubular neighborhoods of incomparable strata are disjoint).

If $i$ is the height of $A$, let $\delta_i = \delta \circ \widetilde{\pi}_{i+1} \circ \cdots \circ \widetilde{\pi}_d$, $d = \dim W$. Then
$$A(\delta_i) = (\widetilde{\pi}_d^{-1} \circ \cdots \circ \widetilde{\pi}_{i+1}^{-1})(A(\delta))$$
since $l_A \circ \widetilde{\pi}_B = l_A$ for $A < B$.
Let $W(\delta)$ be the union of $A(\delta_{\operatorname{ht}A})$ over all the strata $A$ of $W$. Then let
$$\widetilde{\pi} = \widetilde{\pi}_0 \circ \cdots \circ \widetilde{\pi}_d$$
and $\pi = \widetilde{\pi}|_{W(\delta)}$.
The next proposition collects some basic properties of $\pi$.

\begin{prop}\label{prop:pi properties} We have:
\begin{enumerate}
\item For each stratum $A$ on $W$ of height $i$,
$$\textstyle \pi^{-1}(A) = A(\delta_i) - \bigcup_{B > A} B(\delta_{\operatorname{ht}(B)}).$$
\item For each $A$ as in (ii), $\pi$ coincides with $\pi_A$ on $\pi^{-1}(A)$.
\item For $A$ as in (ii) and each $y \in A$, we have:
    $$\pi^{-1}(y) = \pi_A^{-1}(y) \cap W(\delta).$$
\end{enumerate}
\end{prop}
\begin{proof} Let $x$ be in $A(\delta_i)$ such that $x \not\in B(\delta_{\operatorname{ht}(B)})$ for any $B > A$. Then $\widetilde{\pi}(x) = (\widetilde{\pi}_0 \circ \cdots \circ \widetilde{\pi}_i)(x_{i+1})$ for $x_{i+1} = (\widetilde{\pi}_{i+1} \circ \cdots \circ \widetilde{\pi}_d)(x)$. We note that $x_{i+1}$ is in $A(\delta)$ and so $\widetilde{\pi}_i(x)$ is in $A$. Since $\widetilde{\pi}_j$ is the identity on $W$ for $j < i$, we conclude that $\widetilde{\pi}(x)$ is in $A$. Conversely, let $x$ be in $\pi^{-1}(A)$. Then $\widetilde{\pi}(x) = \widetilde{\pi}_A(x_{i+1})$ is in $A$ for $x_{i+1}$ defined as before. But that is possible only if $x_{i+1} \in A(\delta)$ or $x \in A(\delta_i)$.

The proof of (ii) is similar: use the fact that $\pi_A \circ \widetilde{\pi}_B = \pi_A$ for $A < B$. (iii) is also seen in the similar way.
%Hence, if $i$ is the height of $A$ and $x$ in $A(\delta_i)$, then we have $(\widetilde{\pi}_{i+1} \circ \cdots \circ \widetilde{\pi}_d)(x)$ is in $A(\delta)$ and so
%$$(\widetilde{\pi}_i \circ \cdots \circ \widetilde{\pi}_d)(x) = (\pi_A \circ \widetilde{\pi}_{i+1} \circ\cdots \circ \widetilde{\pi}_d)(x) = \cdots = \pi_A(x).$$
%Since $\widetilde{\pi}_j|_W$ is the identity for any $j$, we have that $\widetilde{\pi}$ (thus $\pi$) coincides with $\pi_A$ on $A(\delta_i)$.
\end{proof}

\begin{corollary}\label{W(delta) interior} In the definition of $\widetilde{\pi}$, the function $\delta$ can be chosen so that the interior of $W(\delta)$ is a neighborhood of $W$.
\end{corollary}
\begin{proof} If $x$ is in $\pi^{-1}(A) - A$, then $\delta_i > 0$ near $x$. Thus, if we assume $\delta_i$ decreases slowly compared to $l_A$, then $\pi^{-1}(A)$ contains a point near $A$ unless the point is already in $B(\delta_{\operatorname{ht}(B)})$ for some $B > A$. That is, $W(\delta)$ contains a neighborhood of $A$. Since the stratification on $W$ is locally finite, such an assumption on $\delta$ need to be made again only finitely strata at a time. The assertion is therefore clear.
\end{proof}

After choosing $\delta$ according to the above corollary, we have:

\begin{corollary} $\pi(W(\delta)) = W$.
\end{corollary}
\begin{proof} This follows from the previous corollary as well as (iii) of the proposition.
\end{proof}

In fact, $\pi$ is much more than a surjection.

\begin{prop}\label{pi prop} For $U$ = the interior of $W(\delta)$, we have $\pi : U \to W$ is a strong deformation retract of $W \subset U$; i.e., there is a homotopy $h_t : U \to U$ from $\pi$ to the identity such that $h_t|_W$ is the identity. Moreover, $h_t$ can be taken to be smooth on $U - W$.
\end{prop}
\begin{proof} For a stratum $A$ of $W$, if $x = H_A(x_0, s)$ and $\delta(x) > 0$, then let
$$\widetilde{\pi}_{A, t}(x) = H_A(x_0, (1 - (1 - t) \lambda_A(x))s)$$
and $\widetilde{\pi}_{A, t}(x) = x$ for any other $x$. Note $\widetilde{\pi}_{A, t}(U) \subset U$. Then define $\widetilde{\pi}_t$ as we defined $\widetilde{\pi}$ from $\widetilde{\pi}_A$. Take $h_t : U \to U$ to be the restriction of $\widetilde{\pi}_t$.
\end{proof}

An important consequence of the above that we will need in the next section is the corollary below (as a reader will see, the proof will go through if $\pi$ is weaker than a strong deformation retract. ``strong'' is used in the next section). 

\begin{corollary}[Relative Poincar\'e Lemma]\label{cor:Poincare} Let $\sigma$ be a continuous $k$-form, $k > 0$ (Remark \ref{continous forms}) on some neighborhood $U$ of $W$ such that $\sigma$ is smooth and closed on $U - W$. If $\sigma|_W = 0$, then, after shrinking $U$, there exists a continuous $(k-1)$-form $\psi$ on $U$ that is smooth on $U - W$ with $d \psi = \sigma$ on $U - W$ and $\psi|_W = 0$.
\end{corollary}
\begin{proof} Assume $U$ is small enough that $h_t : U \to U$ for $h_t$ in the preceding proposition. Let $p : U \times [0, 1] \to U$ denote the projection and, since $p$ has compact oriented fibers, we can define the integration along fibers $p_*$ by: for a smooth form $\alpha$ not involving $dt$, $t$ the coordinate on $[0, 1]$, $p_*(\alpha) = 0$ and
$$p_*(f \, dt \wedge \alpha) = \left( \int_0^1 f(\cdot, t) \, dt \right) \alpha.$$
This local formula shows that $p_*$ is in fact defined for any continuous form $\alpha$ on $U \times [0, 1]$ (remark: integrability is enough). From the local formula, by differentiation under the integral sign, we see that for any smooth form $\alpha$ on $(U - W) \times [0, 1]$,
$$d(p_* \alpha) = \alpha|_{U \times 1} - \alpha|_{U \times 0} - p_*(d \alpha).$$
Let $J = p_* \circ h^*$, which is defined on $U$. Then, taking $\alpha = h^* \sigma$ on $U - W$ for $h : U \times [0, 1] \to U, \, (x, t) \mapsto h_t(x)$, we have: on $U - W$,
$$h_1^* \sigma - h_0^* \sigma = J d \sigma + dJ \sigma.$$
Let $\psi = J \sigma$, which is defined on $U$. Now, we have $h_1^* \sigma = \sigma$ since $h_1$ is the identity and we have $h_0^* \sigma = 0$ on $U$ since $h_0|_U = \pi$ and $\sigma|_W = 0$. Since $d \sigma = 0$ on $U - W$ by assumption, we get $\sigma = d J \sigma$ on $U - W$. Finally, we have $\psi|_W = 0$ since $h_t(W) \subset W$.
\end{proof}

\begin{proof}[Proof of Theorem \ref{thm:etale map to a special fiber}]

There is a small $0 \ne t \in \mathbb{C}$ such that $p^{-1}(t)$ lies in the interior of $W(\delta)$. 
Indeed, if this is not the case, we can choose a sequence of $x_i$ not in the interior of $W(\delta)$ such that $p(x_i) \to 0$. 
We have $x_i$ is in the pre-image under $p$ of some closed ball, which is compact by proper-ness. Then, passing to a subsequence, we can assume $x_i$ is convergent with limit in $W = p^{-1}(0)$. But then $x_i$ lies in the interior of $W(\delta)$ for some large $i$, a contradiction.

Choose Whitney stratifications on $\mathfrak{X}$ and $\mathbb{C}$ so that $p$ maps each stratum submersively to a stratum on $\mathbb{C} - 0$; this is possible by \cite[Part 1, Theorem 1.7]{StratifiedMorse}.
 By \cite[Lemma 7.3.]{Mather}, we can assume that $\pi_A$ restricted to a stratum on $\mathfrak{X}$ whose closure intersecting $W$ is a submersion. Making $t$ smaller, we can also assume $p^{-1}(t)$ intersects only strata whose closures intersect $W$.

Fix a stratum $A$ of $W$ of height $i$. We note $\pi^{-1}(A) \cap \mathfrak{X} - W$ has a Whitney stratification induced from $\mathfrak{X}$, as in the proof of Proposition \ref{conical structure}, after modifying $\delta$ if needed (i.e., $\delta$ should be a general choice relative to $\mathfrak{X}$). Let $p_A = p|_{\pi^{-1}(A) \cap \mathfrak{X} - W}$. Then, by Lemma \ref{second isotopy lemma}, there is a neighborhood $U$ of $t$ in $\mathbb{C} - 0$ such that:
$$\xymatrix{p_A^{-1}(t) \times U \ar[rd]_{\pi_A \times \operatorname{id}} \ar[rr]_{\sim} & & p_A^{-1}(U) \ar[ld]^{(\pi_A, p_A)} \\
& A \times U.}\\$$
By Remark \ref{Hironaka}, the horizontal homeomorphism can be assumed to be a stratum-wise diffeomorphism. Then the diagram implies that $\pi_A$ (thus $\pi$) is a submersion on each stratum of $p_A^{-1}(t)$ to $A$ (note, as in the second proof of Lemma~\ref{second isotopy lemma}, $(\pi_A, p_A)$ is a submersion on each stratum.). By \cite[Proposition 7.1.]{Mather}, we have a system of tubular neighborhoods of the strata of $\mathfrak{X} - W$ compatible with $\pi$ where $\pi$ is defined. Then this system restricts to $p^{-1}(t)$. Since, by the above, $\pi : \pi^{-1}(A) \cap p^{-1}(t) \to A$ is a stratum-wise submersion onto each stratum $A$ on $W$, by Thom's first isotopy lemma (to be precise, its abstract version \cite[Corollary 10.2]{Mather}), $\pi|_{p^{-1}(t)}$ is a locally trivial fibration over each stratum; thus is a finite covering over each stratum, as the fibers are finite sets (see (iii) in Proposition \ref{prop:pi properties}). The theorem is then the consequence of the following comparison theorem in \'etale topology.
\end{proof}

\begin{theorem}[SGA 4, Expos\'e XI, Th\'eor\`em 4.3. (iii)]\label{topological coverings of a complex variety} 
Let $X$ be an algebraic variety over $\mathbb{C}$.
Then the functor $X \mapsto X(\mathbb{C})$ determines an equivalence between the category of finite \'etale coverings of $X$ and the category of finite (topological) coverings of $X(\mathbb{C})$ in the classical topology.

Moreover, scheme structures on coverings are induced naturally; i.e., in a manner analogous to manifold structures are induced naturally on coverings of manifolds. In particular, if a covering a priori had a manifold structure and the covering map is smooth, then an induced smooth variety structure would coincide with it.
\end{theorem}
The above theorem is often stated for a smooth complex variety but smoothness is actually not needed.

\begin{proof} The second statement is a direct consequence of the proof of the first statement in \emph{loc.cit.}.
\end{proof}

\section{Moment maps and integrable systems}\label{sec:moment maps}

In this section, we shall give applications of the collapsing map constructed in the previous section, including the construction of an integrable system in the manner of \cite{HK15}.

\subsection{Moment map}
\label{sec:moment map setup}

The notion of a moment map appears in the context of Hamiltonian group actions (Hamiltonian dynamics) in symplectic geometry. A basic case is that of a torus action and in particular, there is a moment map for a toric variety in a projective space (possibly non-normal). By the convexity theorem, the image is a convex polytope, the polytope $\Delta(\mathcal{O}_W(1))$ associated to the toric variety $W$ (Definition~\ref{delta 1}), up to a sign depending on the convention. In fact, we may define a moment map over an algebraically closed field of characteristic zero and we will show that the convexity result remains valid. 

Then composing the moment map a toric variety $W$ with a map $X\to W$ gives a map from the variety $X$ to the polytope of $W$. Over $\mathbb{C}$ at least, the components of the moment map constitute an integral system.
As in \cite{HK15}, we then see that the pull-back of the components along the surjective map $X\to W$ gives an integral system on the original variety $X$.

\begin{definition}\label{delta 1} Let $W$ be a (possibly non-normal) toric variety that appears as the closure of an orbit in a projective space. Precisely, let $\mathbb{G}_m^r$ act on $\mathbb{P}^q$ so that the homogeneous coordinates $x_i$'s have weights $a_i$; i.e., $a_i$ is in $\mathbb{Z}^r$ and $t \cdot x_i = t_1^{a_{i, 1}} \cdots t_r^{a_{i, r}} x_i$ for $t$ in ${k^*}^r = \mathbb{G}_m^r(k)$. Then $W = \overline{\mathbb{G}_m^r \cdot [1 : \cdots : 1]}$ is a (not-necessarily-normal) toric variety.

Then the convex polytope $\Delta(\mathcal{O}_W(1))$ associated to a toric variety $W \subset \mathbb{P}^q$ (or more precisely to $\mathcal{O}(1)$) is defined as the convex hull of $\{ a_0, \dots, a_q \}\subset \mathbb{R}^r$.
\end{definition}

For a moment, assume the base field $k = \mathbb{C}$ is the field of complex numbers. Then the real $r$-torus $(S^1)^r$ also acts on $\mathbb{P}^q$ in a way it preserves the Hermitian structure (i.e. it acts as unitary operators).

Let $\mathfrak{t}_r = \operatorname{Lie}((S^1)^r)$ be the real Lie algebra of the real $r$-torus, which again acts on $\mathbb{P}^q$. Then there is the moment map \cite[Example 3.4.2.]{Moment_GIT}:
\[
\mu : \mathbb{P}^q \to \mathfrak{t}_r^*
\quad \text{given by} \quad 
\langle \mu(x), \xi \rangle = \frac{\operatorname{Im}(\widetilde{x}, \xi \cdot \widetilde{x})}{|\widetilde{x}|^2}
\]
where $\mathfrak{t}_r^*$ is the dual real vector space, $\widetilde{x}$ is a representative in $\mathbb{C}^{q+1}$ of the point $x$ and $(\ ,\ )$ is the standard inner product. Explicitly, 
$$\mu([z_0 : \cdots : z_q]) = \frac{1}{|z|^2} \sum_{j = 0}^q |z_j|^2 a_j$$
since $\mathbb{G}_m^r$ acts on $\mathbb{C}^{q+1}$ with weights $-a_0, \dots, -a_q$. The explicit formula can be used to define a moment map for a base field other than $\mathbb{C}$ if it still has characteristic zero. Indeed, assume the base field $k$ is an algebraically closed field of characteristic zero. Then we can find a real closed field $k_{\mathbb{R}}$ such that $k = k_{\mathbb{R}}(i), i^2 = -1$ (see for example \cite{Nagata} \S 5.1. Exercise 7.) Thus, as a $k_{\mathbb{R}}$-vector space, $k$ is $k_{\mathbb{R}}^2$ with respect to the basis $1, i$. Then we can define the Euclidean norm $|z|$ of $z \in k$ by $|z| = \sqrt{x^2 + y^2}$ when $z = x + i y$. We then define the moment map
$$\mu : \mathbb{P}^q \to \mathfrak{t}_r^*$$
by the above formula, where we used the notation $\mathfrak{t}_r = k_{\mathbb{R}}^r$ and viewed each $a_j$ as a linear functional on $\mathfrak{t}_r$.

\

For an $(S^1)^r$-submanifold of a complex projective space, the next lemma is a part of the Kempf--Ness theorem stating that a GIT quotient can be viewed as a symplectic quotient (cf. \cite[Theorem 8.3]{GIT}). A key part of the theorem is that the moment map gives a characterization of semistable points. The version here is more general in that it holds for an arbitrary algebraically closed field of characteristic zero (we do not know how to extend it to positive characteristic).

We recall that a closed point $x$ in $\mathbb{P}^q$ is said to be \emph{semistable} if either of the following equivalent conditions is met:
\begin{itemize}
\item For a nonzero vector $\widetilde{x} \in k^{q+1}$ representing $x$, the orbit closure $\mathbb{G}_m^r \cdot \widetilde{x}$ does not contain the origin $0$.
\item There exists a nonconstant $\mathbb{G}_m^r$-invariant homogeneous polynomial in the homogeneous coordinates on $\mathbb{P}^q$ that does not vanish at $x$.
\end{itemize}

\begin{lemma}
\label{semistable moment characterization} 
A point $x$ in $\mathbb{P}^q$ is a semistable point if and only if $\mu^{-1}(0)$ intersects the orbit closure $\overline{\mathbb{G}_m^r \cdot x}$.
\end{lemma}
\begin{proof}
First suppose $\mu^{-1}(0)$ intersects $\overline{\mathbb{G}_m^r \cdot x}$; i.e., there is a point $y$ in $\overline{\mathbb{G}_m^r \cdot x}$ such that $\mu(y) = 0$. We shall show $y$ is semistable. Suppose otherwise. Let $v$ be a nonzero vector in $k^{q+1}$ that is a representative of $y$. By the Hilbert--Mumford criterion, we can find a homomorphism $\rho : \mathbb{G}_m \to \mathbb{G}_m^r$ such that $0$ is in the (Zariski-)closure of $\rho(\mathbb{G}_m) v$. Replacing $\rho$ by $\rho^{-1}$ if needed, we can assume $\rho(t) v \to 0$ as $t \to \infty$, in the sense of specialization. Since $\mathbb{G}_m^r$ acts diagonally, the action of $\mathbb{G}_m$ through $\rho$ is also diagonal:
\[
\rho(t)v = \operatorname{diag}(t^{-b_0}, \dots, t^{-b_q})v
\]
for some integers $b_i $. 
Then, for each $i$, $t^{-b_i} v_i \to 0$ as $t \to \infty$. This means $v_i = 0$ if $b_i \le 0$. On the other hand, the moment map corresponding to the $\mathbb{G}_m$-action through $\rho$ is $\rho^* \circ \mu$ where $\rho^* : \mathfrak{t}_r^* \to \mathfrak{t}_1^*$ is the transpose (pullback) of $\rho$ and then $(\rho^* \circ \mu)(y) = |v|^{-2} \sum b_j |v_j|^2 = 0$ since $\mu(y) = 0$.
This is possible only if $v_j = 0$ for all $j$ such that $b_j > 0$. 
Thus, $v = 0$, a contradiction. Hence, $y$ is semistable, proving the claim. Now, if $x$ is not semistable, then $\overline{\mathbb{G}_m^r \cdot x}$ lies in the unstable locus $X^{us}$ which is the complement of the semistable locus and thus $y$ is also in $X^{us}$, a contradiction. Hence, $x$ is semistable.

Next, suppose $x$ is semistable and $y$ a point in $\overline{\mathbb{G}_m^r \cdot x}$ such that $\mathbb{G}_m^r \cdot y$ is closed ($y$ is usually called a polystable point). Let $v$ be a nonzero vector in $k^{q+1}$ representing $y$ and define the function $\psi_v$ on $(k^*)^r$ by $\psi_v(t) = |t \cdot v |^2.$  
Since $0 \not\in \mathbb{G}_m^r \cdot v$, we see that $\psi_v$ is bounded from below by a nonzero number. Since $\mathbb{G}_m^r \cdot v$ is closed, $\psi_v$ admits a minimum $t$ (although $k_{\mathbb{R}}$ may not be Cauchy-complete, a continuous semialgebraic function on a bounded closed semialgebraic set over $k_{\mathbb{R}}$ admits a minimum and a maximum).

Replacing $v$ and $y$ by $t \cdot v$ and $t \cdot y$, without loss of generality, we shall assume $t = 1$, the identity element. 
We shall show $\mu(y) = 0$; for that, we shall show $(\rho^* \circ \mu)(y) = 0$ for every one-parameter subgroup $\rho : \mathbb{G}_m \to \mathbb{G}_m^r$. Recall that $k = k_{\mathbb{R}}(i)$ and then define the function $f : (0, \infty) \to k_{\mathbb{R}}$ by $f(t) = \psi_{v}(\rho(t)) = \sum_{j=0}^q t^{-2b_j} |v_j|^2$ where $b_i$ are as above. 
Then we have: 
$$
\left.\frac{d}{dt} \right|_{t=1} f(t) = -2 \sum_{j=0}^q b_j |v_j|^2 = -2|v|^2 (\rho^* \circ \mu)(y).
$$
Since $f$ has a minimum at $t = 1$, the above has to vanish. That is, $(\rho^* \circ \mu)(y) = 0$, proving the claim.
\end{proof}

Given a semigroup $\mathfrak{S} \subset \mathbb{N} \times \mathbb{Z}^r$, let
\[
\Delta_{\mathbb{Q}}(\mathfrak{S}) = \{ 0 \} \cup \{ u/n : (n, u) \in \mathfrak{S}, n > 0 \} \subset \mathbb{Q}^r.
\]
Then $\Delta_{\mathbb{Q}}(\mathfrak{S})$ is easily seen to be a convex set in $\mathbb{Q}^r$ (see \cite[\S. 2.1.]{algebraic_tex}). Now, we apply this to the \emph{semigroup of $W$} by which we mean this. If $R$ is the homogeneous coordinate ring of $W$, then $R$ has the weight decomposition by the torus action $\mathbb{G}_m^{r+1}$, where the first $\mathbb{G}_m$ corresponds to the $\mathbb{N}_0$-grading on $R$. Then each weight space $R_{n, u}$ has dimension one, since if $f, g$ are in $R_{n, u}$, then $f/g$ is a rational function and that rational function is $\mathbb{G}_m^r$-invariant on an open dense orbit; i.e., a constant function. In other words, $R$ is the semigroup algebra of the semigroup $\mathfrak{S}_W = \{ (n, u) \mid R_{n, u} \ne 0 \}$ (= $\operatorname{sp}(R)$ in the notation of \S \ref{sec:Rees algebra}).

Since $\mathfrak{S}_W$ is finitely generated, we can write $\Delta_{\mathbb{Q}}(\mathfrak{S}) = \{ u \in \mathbb{Q}^r \mid \langle u, v_i \rangle \ge -c_i \}$ for some rational $v_i$'s and $c_i$'s and then define $\Delta(\mathfrak{S}) = \{ u \in k_{\mathbb{R}}^r \mid \langle u, v_i \rangle \ge -c_i \}$. Note $\Delta(\mathfrak{S}_W) = \Delta(\mathcal{O}_W(1))$, since $\mathfrak{S}_W$ is generated by the degrees of weight vectors generating $R$.

%Recall that $\Delta_{\mathbb{Q}}(S)$ was defined in \S \ref{sec:Toric degenerations}.

\begin{prop}[convexity theorem]\label{convexity theorem} Let $W \subset \mathbb{P}^q$ be a possibly non-normal toric variety as above (Definition \ref{delta 1} at the beginning). Then we have $$\mu(W) \cap \mathbb{Q}^r = \Delta_{\mathbb{Q}}(\mathfrak{S})$$
where $\mathfrak{S}$ is the semigroup of $W$. In particular, $\mu(W)$ is convex and coincides with $\Delta(\mathcal{O}_W(1)) = \Delta(\mathfrak{S})$.
\end{prop}
\begin{proof} Let $R$ be the homogeneous coordinate ring of $W$; i.e., the semigroup algebra of $\mathfrak{S}$. For each point $\lambda \in \mathbb{Q}^r$, we have:
\[
\lambda \in \Delta_{\mathbb{Q}}(\mathfrak{S}) \Leftrightarrow (n, n \lambda) \in \mathfrak{S} \text{ for some $n \in \mathbb{N}_{> 0}$} \Leftrightarrow R^{\lambda}_+ := \bigoplus_{n=1}^{\infty} R_{n, n\lambda} \ne 0.
\]
Thus, we need to show $R^{\lambda}_+ \ne 0 \Leftrightarrow \mu^{-1}(\lambda)$ intersects $W$. First suppose $\lambda = 0$. If $R^{\lambda}_+ \ne 0$, then there exists a nonconstant $\mathbb{G}_m^r$-invariant homogeneous polynomial $f$ that does not vanish identically on $W$. In particular, $W$ admits a semistable point $x$. Then, by Lemma \ref{semistable moment characterization}, $\mu^{-1}(0)$ intersects the orbit closure $\overline{\mathbb{G}_m^r \cdot x}$. Since $\overline{\mathbb{G}_m^r \cdot x} \subset W$, that means $\mu^{-1}(0)$ intersects $W$. Conversely, if there is a point $x$ in $\mu^{-1}(0) \cap W$, then $\mu^{-1}(0)$ intersects the orbit at $x$, a fortiori, the orbit closure at $x$; so, again, by Lemma \ref{semistable moment characterization}, $x$ is semistable. This proves the case when $\lambda = 0$.

Next, assume $\lambda$ is in $\mathbb{Z}^r$ and consider the $\mathbb{G}_m$-action on $R$ given by $t \cdot_{\lambda} f = t^{-\lambda} t \cdot f$ for degree-one elements $f$ in $R$. Then $R^{\lambda}_+$ is the $\cdot_{\lambda}$-invariant subspace of $R_+$. With respect to $\cdot_{\lambda}$, the moment map looks $\mu_{\lambda} = \mu - \lambda$. 
It follows that $\mu_{\lambda}^{-1}(0) = \mu^{-1}(\lambda)$ and so the claim follows from the early case. For an integer $n > 0$, we have $n \Delta_{\mathbb{Q}}(\mathfrak{S}) = \Delta_{\mathbb{Q}}(\mathfrak{S}^{[n]})$ where $\mathfrak{S}^{[n]}$ is the $n$-th Veronese of $\mathfrak{S}$ and so the rational case follows from the integral case, which completes the proof of the first assertion. Finally, we have $\mu(W) \subset \Delta(\mathfrak{S})$ from the definition of $\mu$. Thus, $\Delta_{\mathbb{Q}}(\mathfrak{S}) = \mu(W) \cap \mathbb{Q}^n \subset \mu(W) \subset \Delta(\mathfrak{S})$, which implies $\mu(W) = \Delta(\mathfrak{S})$ since $\mu(W)$ is closed.
\end{proof}

\begin{remark}[on convexity] For manifolds, a more general version of the above Proposition \ref{convexity theorem} is due to Atiyah and Guillemin--Sternberg in the abelian case and due to Kirwan in the non-abelian case, see e.g. \cite[Theorem 8.1.5]{Moment_GIT}. For singular $W$, that $\mu(W)$ is convex is a theorem of Atiyah \cite[Theorem 8.2.1]{Moment_GIT}.
\end{remark}

\begin{example} Consider the elliptic curve example. In that case, the generators of the semigroup $\mathfrak{S}$ are $\nu(x) = (1, 1), \, \nu(y) = (1, 0), \, \nu(z) = (1, 3)$ and thus we have:
$$\mathbb{C}^* \to \mathbb{P}^2, \, t \mapsto [t : 1 : t^3].$$
The closure of the image is $W = \Proj(\mathbb{C}[\mathfrak S]) = V(y^2 z = x^3)$.
The weight of the $S^1$-action is $(1, 0, 3) \in \mathbb{Z}^3$ and so 
$$\mu([z_0 : z_1 : z_2]) = \frac{1}{|z|^2}(|z_0|^2 + 3|z_2|^2).$$
Hence, the moment image of $W$ is $[0, 3]$. Also, we have $\varphi : [x : y : z] \mapsto [y^2 z : y^3 : z^3]$ and so
$$(\mu \circ \varphi)([x : y : z]) = \frac{1}{|y^2z|^2 + |y^3|^2 + |z^3|^2}(|y^2z|^2 + 3|z^3|^2).$$
The image is again the interval $[0, 3]$ as expected.
\end{example}

\begin{example} Let $X$ be an elliptic curve over $\mathbb{C}$ and suppose we have a toric degeneration of $X$ to some $W$. Also, suppose we are given a morphism $\varphi : X \to W$. Since $X$ is normal, $\varphi$ factors as $X \overset{f}\to \mathbb{P}^1 \to W$ where the second map $\mathbb{P}^1 \to W$ is the normalization. If $f$ is finite or equivalently $\varphi$ is finite, then, by the Riemann--Hurwitz formula \cite[Ch. IV, Exercise 2.5. (a)]{Hart}, we see that $f$ ramifies over more than two points. Now, the normalization of a curve is unramified (since birational) and so we conclude that $\varphi$ necessarily ramifies over more than two points. In particular, $\varphi$ cannot be unramified over an open dense orbit if $\varphi$ is finite.
\end{example}

%This ends a quick discussion of moment maps on toric varieties. Finally, from either Theorem \ref{thm:etale map to a special fiber} or Corollary \ref{cor:dominant morphism exists}, we have a map $\varphi : X - B \to W$ with (Zariski-)dense image for some closed subset $B$ ($B$ is empty in the first case). Composing it with $\mu$ we get:
%$$\varphi_{\mu} : X - B \to \mathfrak{t}_r^*$$
%whose image is contained in the polytope $\Delta_W(\mathcal{O}(1))$ associated to $W \subset \mathbb{P}^q$. Moreover, if $W$ is the special fiber of a toric degeneration of $X$, then the image is dense in the polytope.

\subsection{Construction of integrable systems}

We now want to recover the construction of an integrable system after Harada and Kaveh \cite{HK15} in our setup. We shall assume the base field $k = \mathbb{C}$ is the complex field.

We first recall some notions from symplectic geometry. 
Let $(M, \omega)$ be a symplectic manifold. To each smooth function $f$ on $M$, by non-degeneracy, we can have a unique vector field $X_f$ such that $i(X_f) \omega = d f$, where $i(X_f)$ denotes the interior product; i.e., $i(X_f)\omega = \omega(X_f, \cdot)$. We can then form a pairing of two smooth functions:
$$\{ f, g \} = \omega(X_f, X_g)$$
called the \emph{Poisson bracket of $f, g$} \cite[Definition 18.5.]{Da_Silva}. 
A calculation in local coordinates shows that it is a Lie bracket that also satisfies the Leibniz rule. We then say smooth functions $f_1, \dots, f_r, r = \frac{1}{2} \dim M$ form a \emph{(completely) integrable system on $M$} if (1) $d f_1, \dots, d f_r$ are linearly independent at each point and (2) they pairwise Poisson-commute; i.e., $\{ f_i, f_j \} = 0$ for any $i, j$ \cite[Definition 18.10.]{Da_Silva}.

We have an extension of the above notion to a stratified space.

\begin{definition}[Integrable system in the stratified sense; cf. \cite{HK15} Definition 1.1.]\label{integrable system in stratified} 
Let $M$ be a manifold and $X \subset M$ a closed subset stratified by submanifolds of $M$. 
A \textit{sympletic form $\omega$ on $X$ in the stratified sense} is a continuous 2-form 
on $M$ such that its restriction to each stratum $A$ of $X$ is symplectic, \emph{i.e.} smooth, closed and nondegenerate.

Then we say real-valued continuous functions $f_1, \dots, f_r$ on $X$ form an \textit{integrable system in the stratified sense with respect to $\omega$} if (1) they are smooth functions on each stratum and (2) on each stratum $A$, they pairwise Poisson-commute with respect to $\omega|_A$ and the span of their differentials at each point has half the dimension of $A$.
\end{definition}

\begin{example} Let $\mu : W \to \mathfrak{t}_r^*$ be the moment map as above. Then the components of $\mu$
$$(f_1, \dots, f_r) : W \overset{\mu}\to \mathfrak{t}^*_r \simeq \mathbb{R}^r$$
form an integrable system on $W$ in the stratified sense, when $W$ is given the stratification by torus orbits.
\end{example}

\begin{remark}[sympletic singular spaces] Our notion of a sympletic structure likely differs from the ones in literature.

There is the notion of a \emph{sympletic variety} due to Beauville \cite[Definition 1.1.]{sympletic_variety}. By defintion, a complex variety is a sympletic variety if it is normal and comes with a closed non-degenerate $2$-form on the regular locus such that $\pi^* \omega$ extends to a holomorphic $2$-form on the entire manifold under any resolution of singularities $\pi : X' \to X$.

The notion of a {\it sympletic stratified space}; i.e., a stratified space with a symplectic structure was introduced by R. Sjamaar and E. Lerman \cite{symplectic_stratified}, who showed that a symplectic reduction with respect to a singular value of the moment map still has a structure of a sympletic stratified space. Among the others, they also showed, in the setup of a sympletic stratified space, that there is a natural Poisson bracket and that Hamiltonian dynamic still works.
\end{remark}

Next, we note that Weinstein's theorem in symplectic geometry continues to hold in the following singular situation, see  \cite[Theorem 3.4.17.]{Secondary_char} or \cite[Theorem 7.4.]{Da_Silva}.

\begin{theorem}[Weinstein]\label{Weinstein} Let $M, W$ be as in the previous section \S \ref{sec:collapsing a family}. Let $\omega_i, i = 0, 1$ be two continuous
%\footnote{By a \emph{continuous $k$-form}, we mean a continuous section of the $k$-th exterior power of $T^* U$. If it is smooth, then it is a differential $k$-form.}
 $2$-forms on some neighborhood $U$ of $W$ such that $\omega_0|_W = \omega_1|_W$. Assume $\omega_0$ is nondegenerate (so $\omega_0$ is a sympletic form if smooth).

If $\omega_i, i = 0, 1$ are locally Lipschitz continuous, as in Lemma \ref{Lipschitz}, and if they are closed smooth forms on $U - W$, then, after shrinking $U$, there exists a homeomorphism $f : U \to U$ such that $f : U - W \to U - W$ is a diffeomorphism, $f : W \to W$ is the identity and $f^* {\omega_1} = \omega_0$ on $U - W$.
\end{theorem}

The original proof uses Moser's trick and goes through without a change due to the relative Poincar\'e lemma (Corollary \ref{cor:Poincare}) holding for $W \subset M$. To be self-contained and for the convenience of the readers, we include the proof following \cite[Proof of Theorem 3.4.17.]{Secondary_char}.

\begin{proof} 
Let $\sigma = \omega_1 - \omega_0$ and $\omega_t = \omega_0 + t \sigma$. For some neighborhood $U$ of $W$, by Corollary \ref{cor:Poincare} (relative Poincar\'e lemma), we can find a continuous one-form $\psi$ on $U$ such that, on $U - W$, $\psi$ is smooth with $d \psi = \sigma$ and $\psi|_W = 0$. Also, by the local formula in the proof of Corollary \ref{cor:Poincare}, we see that $\psi$ is locally Lipschitz since $\sigma$ is (note $h$ there is locally Lipschitz). Since $\omega_t|_W = \omega_0|_W$ is non-degenerate, shrinking $U$, we can have that $\omega_t, 0 \le t < 1$ are all non-degenerate. The non-degeneracy means we can find unique continuous vector fields $\xi_t$ on $U$ such that $i(\xi_t) \omega_t = -\psi$, where $i(v)$ denotes the interior product; i.e., $i(v) \omega_t = \omega_t(v, \cdot)$. Note $\xi_t|_W = 0$ since $\psi|_W = 0$. Also, $\xi_t$ is locally Lipschitz uniformly with respect to $t$ since $\omega_t$ and $\psi$ are locally Lipschitz.

The existence and uniqueness theorem for ordinary differential equations applies to a vector field that is locally Lipschitz uniformly with respect to time \cite[Theorem 2.2.]{teschl}. 
Thus, we get the unique solution $f_t : U \to U$ to the ODE with the initial condition:
$$
\frac{d}{d t}f_t(x) = \xi_t(f_t(x)), \, f_0(x) = x.
$$
We know that $f_t$ is defined for $t$ in some interval around $0$ and is a homeomorphism (as $f_{-t}$ is the inverse). 
Let $\Omega \subset \mathbb{R} \times U$ be the subset consisting of all $(t, x)$ such that $x \in U$ and $t$ in the maximal interval for which the solution $f_t(x)$ exists with the initial condition $f_0(x) = x$. If $\xi_t$ is independent of $t$, then we know $\Omega$ is open \cite[Theorem 6.1.]{teschl}.
But an argument similar to the one showing the openness in that case still carries over to show that $\Omega$ here is also open.

Since $\xi_t|_W = 0$, we have that $\mathbb{R} \times W \subset \Omega$. Thus, shrinking $U$, we can assume $[-1, 1] \times U \subset \Omega$; that is, $f_t : U \to U$ is defined for all $-1 \le t \le 1$. On $U - W$, the right-hand side of the above ODE is smooth and thus, by uniqueness, we also have that $f_t : U - W \to U - W$ is a diffeomorphism.

For any smooth form $\alpha$ on $U - W$, we have:
$$\frac{d}{dt} f_t^* \alpha = f_t^* ((i(\xi_t) d + d i(\xi_t))\alpha).$$
(it comes from a Lie derivative; see \cite[Exercise after Definition 6.3.]{Da_Silva}). Then we have:
\begin{align*}
\frac{d}{dt} f_t^* \omega_t &= \left( \frac{d}{ds} f_s^* \omega_t \right) |_{s = t} + \left( \frac{d}{ds} f_t^* \omega_s \right) |_{s = t} = f_t^*(di(\xi_t) \omega_t) + f_t^* \sigma \\
&= f_t^*(-d \psi + \sigma) = 0.
\end{align*}
That is, on $U - W$, we have $f_t^* \omega_t$ is independent of $t$ and so $f_1^* \omega_1 = \omega_0$. Take $f = f_1$.
\end{proof}

%It is desirable to establish Weinstein's theorem over an algebraically closed field of characteristic zero, which, however, requires the avoidance of the use of a differential equation.

The above Weinstein's theorem combined with the construction in the previous section, we get one of the main results of this section. It is stated in a slightly complicated way; we wanted to emphasize the fact that, for a given fixed $W$, roughly any general $X$ would work.

\begin{prop}\label{integrable system exists proposition} Let $M$ be a smooth complex variety and $W \subset M$ a closed subvariety. Fix 
a complex-algebraic Whitney stratification on $W$ and assume there is a symplectic form $\omega$ on $M$ such that $\omega|_W$ is a symplectic form in the stratified sense (Definition \ref{integrable system in stratified}) and there is a moment map $\mu$ on $W$ with respect to $\omega|_W$. Let $\pi : W(\delta) \to W$ be the retraction constructed in the previous section.

Then there exists a neighborhood $U$ of $W$ in $W(\delta)$ in the classical topology such that for each Whitney stratified locally closed subset $X \subset U - W$ that is in general position relative to $\delta$ (in the sense the proof of Theorem \ref{thm:etale map to a special fiber} goes through), we have that $X$ admits an integrable system with respect to $\omega|_X$ in the sense of Definition \ref{integrable system in stratified} above.
\end{prop}
\begin{proof} We shall say a map $\pi$ is a \emph{stratum-wise surjective submersion} if for each stratum $A$, the restriction to $A$ is a surjective submersion. 

Let $U$ be the interior of $W(\delta)$, $\omega_1 = \pi^* \omega$ and $\omega_0 = \omega$. Then $\omega_1$ is locally Lipschitz since $\pi$ (Lemma \ref{Lipschitz}) and $\omega$ are so. By Theorem \ref{Weinstein} applied to $\omega_1, \omega_0$, after shrinking $U$, we find a diffeomorphism $f : U - W \to U - W$ such that $f^* \pi^* \omega = \omega$ on $U - W$. Then the composition
$$\varphi : X \overset{f}\to U - W \overset{\pi}\to W$$
is a stratum-wise surjective submersion. Thus, the composition $\mu \circ \varphi$ is the claimed integrable system.
\end{proof}

We now obtain the following version of the main result of \cite{HK15}.

\begin{theorem}[Harada--Kaveh type Theorem]\label{thm:improved Harada-Kaveh} Let $M$ be a smooth complex variety (e.g., $M = \mathbb{P}^n \times \mathbb{A}^1$) and $\mathfrak{X} \subset M$ a closed subvariety. Assume we have a surjective morphism $p : \mathfrak{X} \to \mathbb{A}^1$ that is the restriction of a proper morphism $M \to \mathbb{A}^1$. Let $W = p^{-1}(0)$ and assume there is a sympletic form $\omega$ on $M$.

Then, for some Whitney stratification on $\mathfrak X$ and for some $t \ne 0$, there exists a continuous surjective map
$$\varphi : X := p^{-1}(t) \to W$$
such that (1) $\varphi^*(\omega|_W) = \omega|_X$ and (2) for each stratum $A$ on $W$, $\varphi^{-1}(A)$ has a natural reducible-variety structure so that $\varphi : \varphi^{-1}(A) \to A$ is finite \'etale. Note $\varphi^{-1}(A)$ need not be a subscheme.

Moreover, if $W$ comes with a moment map $\mu : W \to \mathfrak{t}^*_r$ (e.g., $\mathfrak{X}$ is a toric degeneration), then
$$(f_1, \dots, f_r) : X \overset{\varphi}\to W \overset{\mu}\to \mathfrak{t}^*_r \simeq \mathbb{R}^r$$
is an integrable system in the stratified sense (Definition \ref{integrable system in stratified}). 
Finally, if $\mathfrak{X}$ is a toric degeneration as in Construction~\ref{locally Grobner toric degeneration}, then the stratification on $W$ may be taken to be an orbit stratification.
\end{theorem}

\begin{proof} The first part is exactly Proposition~\ref{integrable system exists proposition} and Theorem~\ref{thm:etale map to a special fiber}. The very last statement is by Construction \ref{locally Grobner toric degeneration} since in that construction, the strata on $W$ are invariant and such refining still gives a Whitney stratification (cf. Remark \ref{canonical stratification}).
\end{proof}

\subsection{Torus embedding description} Given the result in the previous subsection, it is useful to revisit the question of a characterization of a map to a toric variety, as a complement to \cite{algebraic_tex}.

We write $W \subset \mathbb{P}^N$ for a not-necessarily-normal toric variety embedded in a projective space (we only consider the projective case for simplicity).

\begin{prop} Let $\Delta = \Delta(\mathcal{O}_W(1))$. Then there is an equivalence of categories between the following two.
\begin{enumerate}[label=(\alph*)]
    \item The category of continuous maps $X \to W$ from arbitrary topological spaces $X$ with dense images, where a morphism from $\varphi : X \to W$ to $\varphi : X' \to W$ is a continuous map $X \to X'$ such that $$\xymatrix{X \ar[dd] \ar[rd]^{\varphi} & \\
    & W \\
    X' \ar[ru]_{\varphi'}}$$ commutes.
    \item The category of continuous maps $X \to \Delta$ from arbitrary topological spaces with dense images, where morphisms are defined similar to (a).
\end{enumerate}
\end{prop}
\begin{proof} Let $W'$ denote the normalization of $W$ and then we let $W'(\mathbb{R}_{\ge 0})$ be the set of all torus-equivariant morphisms $\Spec(\mathbb{C}[(\mathbb{R}_{\ge 0}, \text{mult})]) \to W'$. We have $W' \subset \mathbb{P}^{q'}$ for some $q'$ and we identity $W'(\mathbb{R}_{\ge 0})$ with a subset of $\{ [c_0 : \cdots : c_{q'}] \mid c_i \in \mathbb{R}_{\ge 0}, \sum c_i = 1 \} \subset \mathbb{P}^{q'}$. Then we know that the restriction of the moment map $\mu_{\mathbb{P}^{q'}} : \mathbb{P}^{q'} \to \mathbb{R}^r$, namely, 
$$\mu_{\mathbb{P}^{q'}} : W'(\mathbb{R}_{\ge 0}) \to \Delta$$
is a homeomorphism. Let $F$ be the functor from (a) to (b) given by $F(\varphi) = \mu_{\mathbb{P}^q} \circ \varphi$. On the other hand, if we are given $\psi : X \to \Delta$, then we have
$$X \to \Delta \simeq W'(\mathbb{R}_{\ge 0}).$$
Using the fact that $W'(\mathbb{R}_{\ge 0})$ is a quotient of $W'$, we then lift the above map to $X \to W'$. Finally, we let $G(\psi)$ be $X \to W' \to W$. Since $W' \overset{\mu_{\mathbb{P}^{q'}}}\to \mathbb{R}^r$ agrees with $W' \to W \overset{\mu_{\mathbb{P}^q}}\to \mathbb{R}^r$, we have that $(G \circ F)(\varphi)$ and $(F \circ G)(\psi)$ coincide with $\varphi, \psi$ on dense subsets; thus, coincide identically. Hence, $F$ is an equivalence of categories.
\end{proof}

%We have a result analogous to the ones in \S \ref{sec:a map to a toric variety section}.

\begin{prop} \label{prop:map factors through torus}
Let $X$ be a topological space. Then there is a natural one-to-one correspondence between the following two
\begin{enumerate}[label=(\alph*)]
    \item A continuous map $\varphi : X \to W$.
    \item Real-valued continuous functions $f_1, \dots, f_r$ on some open dense subset $U \subset X$ such that $U \to \mathbb{G}_m^r$ defined by $f_i$'s extends to $X \to W$, where $\mathbb{G}_m^r \to W$ is the orbit map onto the open dense orbit.
\end{enumerate}
\end{prop}
\begin{proof} From (a) to (b). Let $O \subset W$ be the open dense orbit. Then any map $\varphi : X \to W$, when restricted to $U = \varphi^{-1}(O)$, factors as
$$\varphi : U \overset{f}\to \mathbb{G}_m^r \simeq O \subset W$$
for some $f$. We can then write $f = (f_1, \dots, f_r)$.
%Then, explicitly, $\varphi|_U$ is given as $$x \mapsto [f^{a_0}(x) : \cdots : f^{a_N}(x)]$$ where we used the multi-index notation and $a_i$ are the weights of the homogeneous coordinates. That is to say, $\varphi|_U$ is given as a torus embedding but in the coordinates $f_i$'s.

From (b) to (a). Clear.
%Like in the usual case, any continuous map $\psi : X \to \mathbb{P}^N$ is defined by giving $N$ sections of $L = \psi^*(\mathcal{O}(1))$ that do not vanish simultaneously (here ``continous'' is with respect to the Zariski topology or the classical topology). Namely, if $s_i$ are such sections, then $\psi$ is defined as:
%$$\psi(x) = [(s_0/s_i)(x) : \cdots : (s_{i-1}/s_i)(x) : 1 : (s_{i+1}/s_i)(x) : \cdots : (s_N/s_i)(x)].$$
%for $x$ in $X_i = \{ s_i \ne 0 \}$.
\end{proof}

%\begin{prop} Let $\varphi : X \to W$ be a map that is continuous with respect to either the Zariski topology or the classical topology and $f_i$ as above and $L = \varphi^*(\mathcal{O}(1))$.

%Then there exists a section $s$ of $L^{\otimes n}$ for some $n > 0$ such that $X_s \subset \varphi^{-1}(O)$ and, after a Veronese embedding, $\varphi$ is given by the sections $s f^{a_0}, \dots, s f^{a_N}$.
%\end{prop}
%\begin{proof} Since $\mathcal{O}_W(1)$ is ample, the open subsets of the form $W_y$, $y$ sections of tensor powers of $\mathcal{O}_W(1)$, generate the Zariski topology of $W$. Since $O$ is open in $W$, we can thus choose a section $y$ of $\mathcal{O}_W(l)$ for some $l > 0$ such that $W_y \subset O$. Take $s = \varphi^*(y)$. Then $\varphi^{-1}(V(y)) = V(s)$ and so $X_s \subset \varphi^{-1}(O)$.
%\end{proof}

%\begin{example} Like before, suppose we are given a toric degeneration of an elliptic curve $X$ to $W$. Note $U = \varphi^{-1}(O)$ is Zariski-open (as the complement is a finite set) and thus $X$ is a projective completion of $U$.
%\end{example}

%\begin{remark} In the case of a normal toric variety, a torus embedding can be intrinsically described by a fan of cones. Now, the question is how to do similar things for maps to toric varieties. For example, if $X$ is a line bundle over $X$, there is already some result (namely, ).
%\end{remark}

\footnotesize{
\bibliographystyle{alpha}
\bibliography{bib}
}

\end{document}